\documentclass[a4paper,12pt, reqno]{amsart}
\usepackage{amssymb,amsmath,latexsym,mathtools,enumitem}
\usepackage{color}
\usepackage{comment}
\usepackage{verbatim}

\usepackage{hyperref}
\hypersetup{	
	colorlinks=true,
	linkcolor=blue,
	citecolor=blue,
}

\allowdisplaybreaks

\def\1{{\bf 1}}

\def\bee{\begin{equation}}
\def\eee{\end{equation}}

  \def\bR {{\mathbb R}}

\def\R {{\mathbb R}}

\def\supp{{\mathop {{\rm supp\, }}}}

\newtheorem{thm}{Theorem}[section]
\newtheorem{lemma}[thm]{Lemma}
\newtheorem{defn}[thm]{Definition}
\newtheorem{prop}[thm]{Proposition}
\newtheorem{corollary}[thm]{Corollary}
\newtheorem{remark}[thm]{Remark}

\numberwithin{equation}{section}

\def\qed{{\hfill $\Box$ \bigskip}}

\def\MM{{\mathcal M}}
\def\BB{{\mathcal B}}

\def\LL{{\mathcal L}}

\def\JJ{{\mathcal J}}

\def\R{{\mathbb R}}

\def\E{{\mathbb E}}
\def\C{{\mathbb C}}

\def\P{{\mathbb P}}
\def\N{{\mathbb N}}

\def\wt{\widetilde}

\def\pff{\noindent{\bf Proof} }

\def\d{\delta_D}
\def\dc{\delta_{D^c}}

\begin{document}

\title[Semilinear equations for non-local operators]{Semilinear equations for non-local operators: Beyond the fractional Laplacian}

\author[Ivan Bio\v{c}i\'c \quad Zoran Vondra\v{c}ek \quad and \quad Vanja Wagner]{Ivan Bio\v{c}i\'c$^*$ \quad Zoran Vondra\v{c}ek$^*$ \quad and \quad Vanja Wagner$^*$}

\thanks{$^*$ Research supported in part by the Croatian Science Foundation under the project 4197.}

\date{}

\begin{abstract}
We study semilinear problems in general bounded open sets for non-local operators with exterior and boundary conditions. The operators are more general than the fractional Laplacian. We also give results in case of bounded $C^{1,1}$ open sets.
\end{abstract}

\maketitle

\bigskip
\noindent {\bf AMS 2020 Mathematics Subject Classification}: Primary 35R11; Secondary 31B25, 31C05, 35J61, 45K05, 60J35

\bigskip\noindent
{\bf Keywords and phrases}: Semilinear differential equations, non-local operators

\smallskip
\section{Introduction}\label{s:intro}
Let $D\subset \R^d$, $d\ge 2$, be a bounded open set, $f:D\times \R\to \R$ a function, $\lambda$ a signed measure on $D^c=\R^d\setminus D$ and $\mu$ a signed measure on $\partial D$. In this paper we study the semilinear problem
\begin{equation}\label{e:4-intro}
	\begin{array}{rcll}
		-Lu(x)&=& f(x,u(x))& \quad \text{in } D\\
		u&=&\lambda& \quad \text{in }D^c\\
		W_Du&=&\mu&\quad \text{on }\partial D.
	\end{array}
\end{equation}

The operator $L$ is a second-order operator of the form $L=-\phi(-\Delta)$ where $\phi:(0,\infty)\to (0,\infty)$ is a complete Bernstein function without drift satisfying certain weak scaling conditions. The operator $L$ can be written as a principal value integral
$$
Lu(x)=\textrm{P.V.}\int_{\R^d}(u(y)-u(x))j(|y-x|)\, dy,
$$
where the singular kernel $j$ is completely determined by the function $\phi$. In case $\phi(t)=t^{\alpha/2}$, $\alpha\in (0,2)$,  $-L$ is the fractional Laplacian $(-\Delta)^{\alpha/2}$ and the kernel $j(|y-x|)$ is proportional to $|y-x|^{-d-\alpha}$.   

The operator $W_D$ is a boundary trace operator first introduced in \cite{BJK} in the case of the fractional Laplacian, and extended to more general non-local operators in \cite{Bio} -- see Subsection \ref{ss:trace} for the precise definition. 

Motivated by the recent preprint \cite{AGCV19} we consider solutions of \eqref{e:4-intro} in the weak dual sense, cf.~Definition \ref{d:wds}, and show that for bounded $C^{1,1}$ open sets this is equivalent to the notion of weak $L^1$ solution as in \cite[Definition 1.3]{Aba15a}.

For the nonlinearity $f$  throughout the paper we assume the condition

\medskip
\noindent 
\textbf{(F)} $f:D\times\R\to\R$ is continuous in the second variable and there exist a function $\rho:D\to[0,\infty)$ and a continuous function $\Lambda:[0,\infty)\to [0,\infty)$ such that $|f(x,t)|\leq \rho(x)\Lambda(|t|)$.
\medskip

Semilinear problems for the Laplacian have been studied for at least 40 years and we refer the reader to the monograph \cite{MV} for a detailed account. The study of semilinear problems for non-local operators is more recent and is mostly focused on the fractional Laplacian, see \cite{FQ12, CFQ, Aba15a, Aba17, BCBF16, BC17, BCBF18, BJK, F20}. One of the important differences between the local and non-local equations is that in the non-local case the boundary blow-up solutions are possible even for linear equations.  To be more precise, there exist non-negative harmonic functions for the operator $L$ that blow up at the boundary. In this paper we will restrict ourselves to the so called moderate blow-up solutions, that is those bounded by harmonic functions with respect to the operator $L$. This restriction is a consequence of the problem \eqref{e:4-intro} itself, namely of the boundary trace requirement on the solution. In this respect we follow \cite{Aba15a,BJK} where the boundary behavior of solutions was also imposed. Note that in \cite{Aba15a} the theory was developed for the fractional Laplacian in a bounded $C^{1,1}$ open set $D$, while \cite{BJK} extends part of the theory to regular open sets. This extension was possible mainly due to potential-theoretic results from \cite{BKK}. 

The goal of this paper is to generalize results from \cite{Aba15a, BJK} and at the same time to provide a unified approach. The first main contribution of the paper is that we replace the fractional Laplacian with a more general non-local operator. This is possible due to potential-theoretic and analytic properties of such operators developed in the last ten years. 
For the most recent development see \cite{BJ20, BL, KKLL, KL20, KL21}. 
Here we single out the construction of the boundary trace operator for the operator $L$ in the recent preprint \cite{Bio}. 
The second main contribution is that we obtain some of the results from  \cite{Aba15a} (which deals with $C^{1,1}$ open sets) for regular open subsets of $\R^d$. To achieve this goal we combine methods from \cite{Aba15a} with those of \cite{BJK}. 

Let us now describe the content of the paper in more detail. In the next section we introduce notions relevant to the paper and recall known results. This includes the notion of the non-local operator $L$, the underlying stochastic process $X=(X_t)_{t\ge 0}$ and its killed version upon exiting an open set, the notion of harmonic function with respect to $X$ (or $L$), and the Green, Poisson and Martin kernel of an arbitrary open subset of $\R^d$. We explain accessible and inaccessible boundary points and its importance to the theory. The boundary trace operator is introduced in Subsection \ref{ss:trace}, cf.~Definition \ref{d:boundary-operator}. The section ends with several auxiliary results about continuity of Green potentials. 

Section \ref{s:semilinear} is central to the paper and contains two main results on the existence of a solution to the semilinear problem \eqref{e:4-intro} in arbitrary bounded open sets.  The first result, Theorem \ref{t:any-f}, can be thought of as a generalization of \cite[Theorem 1.5]{Aba15a}. It assumes the existence of a subsolution and a supersolution to the problem \eqref{e:4-intro} and gives several sufficient conditions for the existence of a solution. As in almost all existence proofs of  semilinear problems, the solution is obtained by using Schauder's fixed point theorem. As a corollary of the third part of that theorem, in Corollary \ref{c:bogdan} we obtain a generalization of the main result of \cite{BJK}. Theorem \ref{t:nonpositive-f} deals with nonpositive function $f$ and is a generalization of \cite[Theorem 1.7]{Aba15a}. The main novelty of our approach is contained in using Lemma \ref{l:approximation} to approximate a nonnegative harmonic function by an increasing sequence of Green potentials. This replaces the approximation used in \cite{Aba15a} which works only in smooth open sets. 

In the last two sections we look at the semilinear problem for $L$ in bounded $C^{1,1}$ open sets and at some related questions. In Section \ref{s:aux-C11} we first recall the notion of the renewal function whose importance comes form the fact that it gives exact decay rate of harmonic functions at the boundary. We then state known sharp two-sided estimates for the Green function, Poisson kernel, Martin kernel and the killing function in terms of the renewal function. Subsection \ref{ss:g-p-potentials} may be of independent interest - there we give the boundary behavior of the Green potential and the Poisson potential of a function of the distance to the boundary. We next provide a sufficient integral condition (in terms of the renewal function) for a function of the distance to the boundary to be in the Kato class. In Subsection \ref{ss:gnd} we invoke a powerful result from \cite{KKLL} to show the existence of generalized normal derivative at the boundary which is used in the equivalent formulation of the weak dual solution. We end the section with a discussion on the relationship of the boundary trace operator $W_D$ with the boundary operator used in \cite{Aba15a,Aba17}.

The last section revisits Theorem \ref{t:nonpositive-f} and Corollary \ref{c:bogdan} in bounded $C^{1,1}$ sets. In case when $f(x,t)=W(\d(x))\Lambda(t)$ for some function $W$,  we give a sufficient and necessary integral condition for (a version of) Theorem \ref{t:nonpositive-f} to hold in terms of $W$, $\Lambda$ and the renewal function. Building on Lemma \ref{l:BJK-126}  we next give a sufficient condition for Corollary \ref{c:bogdan} to hold in a bounded $C^{1,1}$ set.
Finally, we end by establishing Theorem \ref{thm:superlinear} that extends  Corollary \ref{c:bogdan} for nonnegative nonlinearities $f$. This result generalizes \cite[Theorem 1.9]{Aba15a}.

The Appendix has two parts. In the first part we provide a proof of Lemma \ref{l:approximation} in a more general context. In the second part, we give quite technical proofs of Propositions \ref{p:green-potential-estimate} and \ref{p:poisson-potential-estimate}. The proof of Proposition \ref{p:green-potential-estimate} is modeled after the proof of \cite[Theorem 3.4]{AGCV19}, while the proof of Proposition \ref{p:poisson-potential-estimate} is somewhat simpler. 

We end this introduction with a few words about notation. Let $D\subset \R^d$ be an open set. Then $C_b(D)$ denotes the family of all bounded continuous real valued functions on $D$, $C_0(D)$ the family of all continuous functions vanishing at infinity (i.e.~$f\in C_0(D)$ if  for every $\epsilon>0$ there exists a compact subset $K\subset D$ such that $|f(x)|<\epsilon$ for all $x\in D\setminus K$),
$C_c^{\infty}(D)$ the family of all  infinitely differentiable functions with compact support,  $\BB(D)$ Borel measurable functions on $D$, and $\BB_b(D)$ bounded Borel measurable functions on $D$. If $\mu$ is a measure on $D$, then $L^1(D,\mu)$ denotes integrable functions, $L^1_{loc}(D,\mu)$ locally integrable functions and $L^{\infty}(D, \mu)$ essentially bounded functions on $D$.   
In case when $\mu$ is the Lebesgue measure on $D$, we simply write $L^1(D)$, $L^1_{loc}(D)$ and $L^{\infty}(D)$.
Denote by $\partial D$ the boundary of $D$,  $\delta_D(x)=\mathrm{dist}(x, \partial D)$ if  $x\in D$, and $\delta_{D^c}(z)=\mathrm{dist}(z, \partial D)$ if $z\in \overline{D}^c$. For $U\subset D$ open,  $U\subset \subset D$  denotes that the closure $\overline{U}$ is contained in $D$. For $A\subset \R^d$, 
 $\MM(A)$ denotes $\sigma$-finite signed measures on $A$ and $|\lambda |$ denotes the variation of $\lambda\in \MM(A)$. For two positive functions $f$ and $g$, $f\preceq g$ means that the quotient $f/g$ stays bounded from above by a positive constant, and $f\asymp g$ that the quotient $f/g$ stays bounded between two positive constants. Finally, unimportant constants in the paper will be denoted by small letters $c$, $c_1$, $c_2$, $\dots$, and their labeling starts anew in each new statement. More important constants we denote by a big letter $C$, where e.g. $C(a,b)$ means that the constant $C$ depends only on parameters $a$ and $b$.
 
\section{Preliminaries}\label{s:prelim}
\subsection{The process and the jumping kernel}\label{ss:process}
Let $X=(X_t, \P_x)$ be a pure jump L\'evy process in $\R^d$, $d\ge 2$, with the characteristic exponent $\Psi:\R^d\to \C$ given by
$$
\Psi(\xi)=\int_{\R^d}\left(1-e^{i \xi \cdot y}+i \xi \cdot y\1_{\{|y|\le 1\}}\right) \nu(dy),
$$
where $\nu$ is a measure on $\R^d\setminus \{0\}$ satisfying $\int_{\R^d}(1\wedge |y|^2)\nu(dy)<\infty$ -- the L\'evy measure .
Thus the Fourier transform of the distribution of $X_t$ is given by
$$
\E_0 e^{i\xi \cdot X_t}= e^{-t\Psi(\xi)}, \qquad \xi\in \R^d,\ \  t>0.
$$
We further assume that $\Psi(\xi)=\phi(|\xi|^2)$ where $\phi:[0,\infty)\to [0,\infty)$ is a complete Bernstein function, cf.~\cite[Chapter 6]{SSV}. This means that 
$$
\phi(\lambda)=\int_0^{\infty}(1-e^{-\lambda t})\mu(t)dt,\quad \lambda>0,
$$
where $\mu:(0,\infty)\to (0,\infty)$ is a completely monotone function such that $\int_0^{\infty}(1\wedge t)\mu(t)dt<\infty$. Thus, in fact, the process $X$ is a subordinate Brownian motion with the L\'evy measure $\nu(dx)=j(|x|)dx$ where $j:(0,\infty)\to (0,\infty)$ is given by 
\begin{equation}\label{e:levy-density}
j(r)= \int_0^\infty (4\pi t)^{-d/2} e^{-r^2/(4t)}\mu(t)dt,\quad r>0.
\end{equation}
We will refer to the function $j$ as the L\'evy density, or the jumping kernel, or simply, the kernel.
The function $j$ is strictly positive, continuous, decreasing and satisfies $\lim_{r\to \infty}j(r)=0$. 

The main example of the process satisfying the assumptions presented in this subsection is the isotropic $\alpha$-stable process in $\R^d$, $\alpha\in (0,2)$. In this case $\Psi(|\xi|)=|\xi|^{\alpha}$, $\phi(\lambda)=\lambda^{\alpha/2}$, and $j(r)=C(d,\alpha) r^{-d-\alpha}$ for some explicit constant $C(d,\alpha)>0$. The isotropic stable process enjoys the exact scaling property which in terms of the complete Bernstein function $\phi(\lambda)=\lambda^{\alpha/2}$  reads  as $\phi(t)/\phi(s)=(t/s)^{\alpha/2}$. A similar property is also needed for the subordinate Brownian motion $X$. Thus we introduce the following weak scaling hypothesis: 
\medskip

\noindent \textbf{(H)}: There exist $R_0>0$, $0<\delta_1\le \delta_2 <1$ and constants $a_1,a_2>0$ such that
\begin{equation}\label{e:wsc}
a_1\left(\frac{t}{s}\right)^{\delta_1} \le \frac{\phi(t)}{\phi(s)}\le a_2\left(\frac{t}{s}\right)^{\delta_2}, \qquad t\ge s\ge R_0.
\end{equation}
The number $R_0$ above is not important: If \textbf{(H)} holds with some $R_0>0$, then it holds with any $R>0$, but with different constants $a_1,a_2$ ($\delta_1$ and $\delta_2$ of course remain the same). 

It is well known that under the assumption \textbf{(H)} the kernel $j$ enjoys sharp two-sided estimates for small $r>0$: For every $R>0$ there exists $C=C(R)\ge 1$ such that
\begin{equation}\label{e:estimate-j}
C^{-1} \phi(r^{-2})r^{-d}\le j(r)\le C \phi(r^{-2})r^{-d}\, , \qquad 0<r<R,
\end{equation}
see for example \cite[(15), Corollary 22]{BGR}.l Moreover, the following properties of $j$ are known: there exists $C=C(\phi)>0$ such that
\begin{equation}\label{e:j1}
j(r)\le C j(r+1), \qquad r>1,
\end{equation}
for every $M>0$ there exists $C=C(M, \phi)>0$ such that
\begin{equation}\label{e:j2}
j(r)\le C j(2r), \quad r\in (0,M),
\end{equation}
cf.~\cite[(2.11), (2.12)]{KSV12}, and there exists $C=C(\phi)>0$ such that
\begin{equation}\label{e:j3}
\left| \left(\frac{d}{dr}\right)^n j(r)\right|\le C j(r), \qquad r\ge 1, \ \ n=1,2,
\end{equation}
cf.~\cite[Proposition 7.2]{BGPR}.
Further, by \cite[Lemma 4.3]{KSV16}, for every $r_0\in (0,1)$, 
\begin{equation}\label{e:j4}
\lim_{\delta\to 0} \sup_{r>r_0}\frac{j(r)}{j(r+\delta)}= 1 \, .
\end{equation}
Properties \eqref{e:j1}--\eqref{e:j4} are used in some of the results that we quote later.

\subsection{The semigroup,  the operator and the potential kernel}\label{ss:semigroup}
For a bounded or nonnegative function $u\in\BB(\R^d)$ and $t\ge 0$, define $P_t u(x):=\E_x [u(X_t)]$. Then $(P_t)_{t\ge 0}$ is the semigroup corresponding to $X$. It is well known that this semigroup has the Feller property, i.e., $P_t:C_0(\R^d)\to C_0(\R^d)$. 

The space $C_c^{\infty}(\R^d)$ of infinitely differentiable functions with compact support is contained in the domain of the infinitesimal generator of the semigroup, and for $u\in C_c^{\infty}(\R^d)$ it holds that
\begin{eqnarray}
Lu(x)&=&\int_{\R^d}\left(u(y)-u(x)-  \nabla u(x) \cdot (y-x)\1_{\{|y-x|\le 1\}}\right)j(|y-x|)\, dy \label{e:L1}\\
&=&\lim_{\epsilon \to 0} \int_{|y-x|>\epsilon} \left(u(y)-u(x)\right)j(|y-x|)\, dy. \label{e:L2}
\end{eqnarray}
In the familiar case of the isotropic stable process the operator $L$ is the fractional Laplacian.

Under our assumption, the process $X$ is also strongly Feller, i.e., $P_t: \BB_b(\R^d)\to C_b(\R^d)$. Indeed, by using \textbf{(H)}, it easily follows that 
$$
\int_{\R^d}\left|\E_0\left[e^{i \xi \cdot X_t}\right]\right| |\xi|^n \, d\xi <\infty
$$
for all $n\ge 1$, cf.~\cite[(3.5)]{KSV18b}. It follows from \cite[Proposition 2.5(xii) and Proposition 28.1]{Sat} that $X_t$ has a density
$$
p(t,x)= (2\pi)^{-d}\int_{\R^d} \cos(x\cdot  \xi) e^{-t\Psi(\xi)}\, d\xi,
$$
which is infinitely differentiable in $x$. This immediately implies the strong Feller property. For $t>0$ and $x,y\in \R^d$ let $p(t,x,y):=p(t, x-y)$. Then $p(t,x,y)$ are  transition densities of $X$ (or the heat kernel) in the sense that $P_t f(x)=\int_{\R^d}p(t,x,y)f(y)\, dy$. 

The process $X$ is transient if it satisfies the Chung-Fuchs condition
$$
\int_0^1 \frac{\lambda^{d/2-1}}{\phi(\lambda)}\, d\lambda<\infty\, .
$$
This condition is satisfied for $d\ge 3$ and we always impose it in case $d=2$.

Under transience one can define the potential kernel (or the Green function) by
$$
G(x)=\int_0^{\infty}p(t, x)\, dt <\infty.
$$
Moreover, under the assumption \textbf{(H)}, one has the following comparability for small $x$, cf.~\cite[Lemma 3.2(b)]{KSV14}: For every $R>0$, there exists $C=C(R)> 1$ such that
\begin{equation}\label{e:estimate-G}
C^{-1}\phi(|x|^{-2})^{-1}|x|^{-d}\le  G(x) \le  C\phi(|x|^{-2})^{-1}|x|^{-d},\qquad |x|\le R.
\end{equation}

\subsection{Harmonic functions}\label{ss:harmonic}
Let $\LL^1=L^1(\R^d, (1\wedge j(|x|))dx)$. For an open set $U\subset \R^d$, let $\tau_U=\inf\{t>0: X_t\notin U\}$ be the first exit time from $U$. 
A function $u\in \LL^1$ is said to be harmonic in an open set $D\subset \R^d$ if for every open $U\subset \overline{U}\subset D$,
$$
u(x)=\E_x\big[u(X_{\tau_U})\big]\, ,\qquad x\in U.
$$
The function $u$ is regular harmonic in $D$ if the above equality holds with $D$ instead of $U$. If $u$ is harmonic in $D$ and $u=0$ in $\overline{D}^c$, then $u$ is said to be singular harmonic.

We say that the scale invariant Harnack inequality is valid if there exists $r_0>0$ and a constant $c=c(r_0)>0$ such that for every $x_0\in \R^d$, every $r\in (0,r_0)$ and every function $u:\R^d\to [0,\infty)$ which is harmonic in the ball $B(x_0,r)$ it holds that
$$
u(x)\le c u(y)\, , \qquad x,y\in B(x_0,r/2).
$$
It is well known that the scale invariant Harnack inequality is valid under the weak scaling condition \textbf{(H)}, cf.~\cite[Theorem 1, Theorem 7]{Grz}. Moreover, nonnegative harmonic functions are locally H\"older continuous, \cite[Theorem 2, Theorem 7]{Grz}. Under condition \eqref{e:j3} it is shown in \cite[Theorem 4.9]{BGPR} that if $u$ is harmonic in an open set $D$, then $u\in C^2(D)$.

For $u\in \LL^1$ define the distribution $\wt{L}u$ by
$$
\langle \wt{L}u, \varphi\rangle :=\int_{\R^d}u(x) L\varphi(x) dx\, , \qquad \varphi\in C_c^{\infty}(\R^d).
$$
Let $D\subset \R^d$ be open. Then $u$ is harmonic in $D$ if and only if $\wt{L}u=0$ in $D$ (as a distribution), cf.~\cite[Lemma 3.1 and Lemma 3.3]{GKL} and \cite[Theorem 3.14 and Theorem 3.16]{Bio}.

\subsection{Transition density and Green function for the killed process}\label{ss:killed}
Let $D\subset \R^d$ be an open set, and $\tau_D=\inf\{t>0:\, X_t\notin D\}$. The killed process $X^D$ (or part of the process $X$) is defined by $X^D_t=X_t$ if $t<\tau_D$, and $X^D_t=\partial$ if $t\ge \tau_D$. Here $\partial $ is an extra point called the cemetery. Every Borel function $f$ on $D$ is extended to $\partial$ by letting $f(\partial)=0$. The semigroup $(P^D_t)_{t\ge 0}$ of the killed process $X^D$ is defined by
$$
P^D_t f(x)=\E_x[f(X^D_t)]=\E_x[f(X_t), t<\tau_D], \quad f\in \BB_b(D).
$$

For $t>0$ and $x,y\in D$ let
$$
p^D(t,x,y)=p(t,x,y)-\E_x[p(t-\tau_D, X_{\tau_D}, y), \tau_D <t].
$$
It is well known that $p^D(t,\cdot, \cdot)$ is symmetric on $D\times D$. By the strong Markov property, $p^D(t,x,y)$ is the transition density of $X^D$, i.e., $P^D_t f(x)=\int_D p^D(t,x,y)f(y)dy$. Moreover, by continuity of $p(t,x,y)$, the Feller and the strong Feller property of $X$, one can show that $p^D(t,x,y)$ is jointly continuous, see \cite[pp. 34-35]{CZ} and \cite[Lemma 2.2 and Proposition 2.3]{KSV19}. Continuity of $p^D(t,x,y)$ implies that the semigroup $(P^D_t)_{t\ge 0}$ is strongly Feller. 

Let
$$
G_D f(x)=\int_0^{\infty}P^D_t f(x)dt=\E_x\int_0^{\tau_D} f(X_t)dt\, .
$$ 
be the potential operator of the killed process $X^D$. This operator admits the symmetric density
$$
G_D(x,y)=\int_0^{\infty}p^D(t,x,y)\, dt\, , \qquad x,y\in D,
$$
which we call the Green function of $X^D$. That is, $G_D f(x)=\int_D G_D(x,y)f(y)\, dy$. We extend the definition of the Green function by $G_D(x,y)=0$ if $x\in D^c$ or $y\in D^c$.  By using Hunt's switching identity, it is standard to derive that for every open $U\subset D$,
$$
G_D(x,y)=\E_x[G_D(X_{\tau_U}, y)], \quad x\in D, y\in D\setminus U.
$$
In particular, for a fixed $y\in D\setminus U$, the function $x\mapsto G_D(x, y)$ is regular harmonic in $U$ and harmonic in $D\setminus \{y\}$. Since harmonic functions are continuous, we get that  $x\mapsto G_D(x, y)$ is continuous in $D\setminus \{y\}$. By symmetry, $y\mapsto G_D(x,y)$ is continuous in $D\setminus \{x\}$.

If $f:D\to[0,\infty]$ such that $G_Df(x)<\infty$, for some $x\in D$, it was shown in  \cite[Remark 2.4]{Bio}  that $G_Df<\infty$ a.e. and $G_Df\in \mathcal{L}^1\cap L^1_{loc} (D)$. In particular, when $D$ is bounded $G_Df\in L^1(D)$.

\subsection{Martin kernel and Poisson kernel}\label{ss:poisson-martin}
From now on we assume that $D$ is a bounded open subset of $\R^d$. The Poisson kernel of $D$ with respect to the process $X$ is defined by
\begin{equation}\label{e:poisson-kernel}
P_D(x,z):=\int_D G_D(x,y)j(|y-z|)\, dy\, \quad x\in D, z\in D^c.
\end{equation}
It is well known and follows from the L\'evy system formula (see \cite[(1.1)]{KSV16}) that $P_D(x,\cdot)$ is the density (with respect to the Lebesgue measure) of the exit distribution of $X$ from $D$ (restricted to $\overline{D}^c$):
\begin{equation}\label{eq:exit_distribution}
\P_x(X_{\tau_D}\in A)=\int_A P_D(x,z)dz\, , \quad A\subset \overline{D}^c.
\end{equation}
Furthermore, it was shown in  \cite[Proposition 3.1]{Bio}  that $P_D(x,z)$ is jointly continuous in $D\times \overline{D}^c$. It is well known that  if $D$ has  a  Lipschitz boundary,  then $\P_x(X_{\tau_D}\in \partial D)=0$, see  \cite[(5.5)]{KSV12},  and thus the equality \eqref{eq:exit_distribution} holds for every $A\subset D^c$. 

We say that $z\in \partial D$ is accessible from $D$ with respect to $X$ if $P_D(x,z)=\infty$ for some (equivalently, every) $x\in D$, and inaccessible otherwise.
The notion of accessible boundary point  was  introduced in \cite{BKK} in the context of the fractional Laplacian. In the very general setting, accessible points were studied in \cite{KSV18a} and \cite{KJ}. It is shown in \cite[Subsection 4.1]{KSV18a} that the subordinate Brownian motion $X$ of our paper satisfies all the assumptions of \cite{KSV18a}, so we are free to use results of that paper. We mention that in case of sufficiently regular boundary $\partial D$ (Lipschitz boundary is fine), all boundary points are accessible.

Let $\partial_M D\subset \partial D$ denote the set of all accessible boundary points of $D$. Fix $x_0\in D$. It is shown in \cite[Lemma 3.4, Theorem 1.1]{KSV18a}  that for every accessible point $z\in \partial D$ there exists 
\begin{equation}\label{e:martin-kernel}
M_D(x,z):=\lim_{y\to z, y\in D} \frac{G_D(x,y)}{G_D(x_0,y)},
\end{equation}
and $x\mapsto M_D(x,z)$ is harmonic with respect to $X^D$ (i.e. singular harmonic with respect to $X$). In fact, the above limit exists for all boundary points $z$, but $x\mapsto M_D(x,z)$ is not harmonic in case of an inaccessible point $z\in\partial D$. The function $M_D(x,z)$ is called the Martin kernel of $D$ with respect to $X$. It is shown in  \cite[Proposition 5.11]{Bio}  (cf.~\cite{BKK} for the case of the fractional Laplacian) that $u:D\to [0,\infty)$ is harmonic with respect to $X^D$ if and only if there exists a nonnegative finite measure $\mu$ on $\partial_M D$ such that
$$
u(x)=\int_{\partial_M D} M_D(x,z)\mu(dz)\, .
$$
In that case $\mu$ is unique. We will use the notation $M_D \mu(x)=\int_{\partial_M D}M_D(x,z)\mu(dz)$. Since $M_D \mu$ is a singular harmonic function with respect to $X$, we have that $M_D\mu\in C^2(D)$, and also by \cite[Remark 5.12]{Bio}, it is in $L^1(D)$.
We note further that $M_D\mu\equiv \infty$ in $D$ if and only if $\mu$ is an infinite measure, see \cite[Corollary 5.13]{Bio}.

For a nonnegative measure $\lambda$ on $D^c$, we define
$$
P_D\lambda (x):=\int_{D^c}P_D(x,z)\lambda (dz),\quad x\in D.
$$
If $\lambda$ is a signed measure on $D^c$ such that $P_D|\lambda|<\infty$  in $D$, then $P_D \lambda$ is defined by the same formula. Note that if $P_D|\lambda|(x)<\infty$, for some $x\in D$, \cite[Corollary 3.11 and Remark 3.6]{Bio}  yield  that $P_D\lambda$ is finite and continuous on the whole $D$, and $P_D\lambda\in L^1(D)$.

We can say something more about the measure that satisfies $P_D|\lambda|<\infty$. Since $P_D(x,z)=\infty$ for $z\in \partial_M D$, $P_D|\lambda|<\infty$ implies that the measure $\lambda$ has no mass on $\partial_MD$ so $\lambda$ can be viewed as a measure on $\R^d\setminus (D\cup \partial_M D)$. Also, $\lambda$ is finite on compact subsets of $\overline D^c$ since for a compact $K\subset\overline D^c$ we have that $K\ni y\mapsto P_D(x,y)$ is bounded and strictly positive.

Let $\nu\in \MM(D)$ and set $u(y)=G_D \nu(y):=\int_D G_D(y,v)\nu(dv)$. Then for $z\in D^c$ we have $u(z)=0$, hence
\begin{eqnarray}
Lu(z)&=&\lim_{\epsilon \to 0}\int_{|y-z|>\epsilon}(u(y)-u(z))j(|y-z|)\, dy=\int_{\R^d}G_D \nu(y)j(|y-z|)\, dy\nonumber\\
&=&\int_D \left(\int_D G_D(y,v)\nu(dv)\right)j(|y-z|)\, dy\nonumber\\
&=&\int_D \left(\int_D G_D(y,v)j(|y-z|)\, dy\right)\, \nu(dv)\nonumber\\ 
&=&\int_D P_D(v, z)\, \nu(dv), \nonumber
\end{eqnarray}
if the last integral absolutely converges.
In particular, if $\nu=\delta_x$ for $x\in D$, where $\delta_x$ is the Dirac measure at $x$, then $u(y)=G_D(x,y)$ and 
$$
L G_D(x, \cdot)(z)=P_D(x,z),
$$
which gives an alternative expression for the Poisson kernel. Further, let $\psi:D\to \R$ be bounded, $u=G_D \psi$ and $\lambda \in \MM(D^c )$. Then
\begin{eqnarray}
\lefteqn{\int_{D^c}Lu(z)\, \lambda(dz)=\int_{D^c}  \left(\int_D P_D(y,z)\psi(y)\, dy\right)\lambda(dz)}  \nonumber \\
&=&\int_D \psi(y)\left(\int_{D^c} P_D(y,z)\, \lambda(dz)\right) dy=\int_D \psi(y) P_D \lambda(y)\, dy. \label{e:Lu}
\end{eqnarray}

\subsection{Boundary trace operator}\label{ss:trace}
Recall that $x_0\in D$ is fixed. Let $u:D\to[-\infty,\infty]$ and let $U\subset\subset D$ be an open Lipschitz set such that $x_0\in U$. Let $\eta_Uu$ be a measure on $\R^d$ defined by
\begin{align*}
	\eta_Uu(A)=\int_A G_U(x_0,z)\left(\int_{D\setminus U} j(|z-y|)u(y)dy\right)dz,\quad A\subset \R^d.
\end{align*}
Following the approach in \cite{BJK} we define the boundary trace operator $W_D$.
\begin{defn}\label{d:boundary-operator}
	If the measures $\eta_U|u|$ are bounded as $U\uparrow D$ and $\eta_Uu$ weakly converge to a measure $\mu$ as $U\uparrow D$, then we denote $\mu$ by $W_Du$, i.e. $W_Du\coloneqq\lim\limits_{U\uparrow D} \eta_Uu$.
\end{defn}
It was shown in \cite[Lemma 5.2]{Bio}  that the measure $W_D u$ is concentrated on $\partial D$. Further, if $\mu$ is a finite signed measure on $\partial_M D$, $\lambda $ a $\sigma$-finite measure on $D^c$ such that $P_D|\lambda|<\infty$, and $f:D\to \R$ such that $G_D|f|<\infty$, then by \cite[Proposition 5.4 and Proposition 5.11]{Bio}, it holds that

\begin{equation}\label{e:W_D(G_D)}
	W_D(M_D\mu)=\mu,\enskip W_D(P_D\lambda)=W_D(G_D f)=0. 	
\end{equation}

\subsection{Some auxiliary results about Green potentials}\label{ss:aux}
We keep assuming that $D$ is a bounded open subset of $\R^d$. Recall that a function $q:D\to [-\infty, \infty]$ is said to be in the Kato class $\JJ$ with respect to $X$ if the family of functions $\{G_D(x,y)|q|(y):\, x\in D\}$ is uniformly integrable (with respect to the Lebesgue measure on $D$). Obviously, if $|v|\le|q|$ and $q\in\JJ$ then $v\in\JJ$. 

Next, we show that a function $q:D\to [-\infty, \infty]$ satisfying
\begin{equation}\label{eq:kato}
\lim_{\epsilon \to 0}\sup_{x\in \R^d} \int_{|x-y|<\epsilon} |q(y)|\phi(|x-y|^{-2})^{-1}|x-y|^{-d}\, dy =0
\end{equation} 
is in the Kato class $\JJ$. Extend the function $q$ to all of $\R^d$ by setting $q(y)=0$ for $y\in D^c$. Since $G_D(x,y)\le G(x,y)$, to show that $q\in\JJ$ it suffices to show that the family of functions $\{G(x,y)|q(y)|:\, x\in \R^d\}$ is uniformly integrable. By using \eqref{e:estimate-G}, one can check that \cite[(24), Lemma 5]{Zhao} holds true. Hence, we can apply \cite[Theorem 1]{Zhao}, which together with \eqref{e:estimate-G} implies that \eqref{eq:kato} is equivalent to
\begin{equation}\label{eq:kato2}
\lim_{t\downarrow0}\sup_{x\in\R^d}\mathbb E_x\left[\int_0^t q(X_s)\, ds\right]=0,
\end{equation}
i.e.~ $q$ is in the classical Kato class $K(X)$ from \cite{CS02} and \cite{C02}. By \eqref{eq:kato}, $q\in L^1(D)$ and therefore \cite[Theorem 2.1(ii)]{C02} implies that $q\in K_{\infty}(X)$, i.e. 
\begin{equation}
\forall\varepsilon>0\ \exists \delta>0\  \forall B\in\mathcal B(\R^d)\text{ such that } \lambda(B)<\delta\ \Rightarrow\ \sup_{x\in\R^d}\int_B |q(y)|G(x,y) dy<\varepsilon.
\end{equation}
cf.~\cite[Definition 2.1(ii)]{CS02}. Furthermore, by \cite[Proposition 2.1]{CS02}, $q\in K_\infty(X)$ implies that $q$ is Green bounded. Together with boundedness of $D$, \cite[Theorem 16.8(iii)]{Sch} gives that the family $\{G(x,y)|q(y)|:\, x\in \R^d\}$ is uniformly integrable, and therefore $q\in\JJ$.

Note that under \textbf{(H)}, the condition \eqref{eq:kato} is satisfied for $q\in \mathcal B_b(\R^d)$, so every bounded function $q$ is in the Kato class $\JJ$.

Recall that the boundary point $z\in \partial D$ is said to be regular (for $D^c$) if $\P_z(\tau_D=0)=1$. The set $D$ is regular if every boundary point is regular. The same proof as in \cite[Proposition 1.31]{BJK} shows that if $D$ is regular, then $q\in \JJ$ if and only if $G_D |q|\in C_0(D)$, and then $G_D q  \in C_0(D)$.

Let $z\in \partial D$ be regular. Then for all $x\in D$, 
$$
\lim_{y\to z, y\in D}G_D(x,y)=0.
$$
A proof of this well-known result can be found in \cite[Proposition 6.2]{KSV17}. The next result is also known -- we include the proof for the sake of completeness.

\begin{lemma}\label{l:continuity-of G1}
Let $D$ be a bounded open subset of $\R^d$. Then $G_D \1\in C(D)$ and $\lim_{x\to z}G_D \1(x)=0$ for every regular boundary point $z\in \partial D$. 
\end{lemma}
\proof Let $(x_n)_{n\in \N}$ be any sequence of points in $D$. Since the constant function \1 is in $\JJ$, the family $\{G_D(x_n, \cdot):\, n\in \N\}$ is uniformly integrable. If $x_n\to x\in D$, then $\lim_{n\to \infty}G_D(x_n, y)=G_D(x,y)$ for a.e.~$y\in D$, hence by Vitali's theorem, see \cite[Theorem 16.6 $(ii)\iff (iii)$]{Sch}, it follows that
$$
\lim_{n\to \infty}\int_D G_D(x_n,y)dy=\int_D G_D(x,y)dy\, ,
$$
proving that $G_D \1 \in C(D)$. If $x_n\to z\in \partial D$ with $z$ regular, then $\lim_{n\to \infty}G_D(x_n,y)=0$ for all $y\in D$. Again by Vitali's theorem we get that $\lim_{n\to \infty}\int_D G_D(x_n, y)\, dy=0$.
\qed

Denote by  $D^{\mathrm{reg}}$ the set of all regular boundary points of $D$. For $\delta>0$, let $D_{\delta}:=\{x\in D:\, \mathrm{dist}(x,\partial D)>\delta\}$.

\begin{lemma}\label{l: lemma_nakon}
	Let $v:D\to[0,\infty)$ be a locally bounded function and  $\rho:D\to[0,\infty)$ such that $G_D\rho\in C(D)$ and $\rho vG_D \1\in L^1(D)$. Then, for every $x\in D$ it follows that
	\begin{align*}
		\lim_{ w\to x}\int_D |G_D(x,y)-G_D(w,y)|\rho(y)v(y)dy=0.
	\end{align*}
\end{lemma}
\proof
Let $r>0$ such that $\overline{B(x,r)}\subset D$ and take a sequence $(x_n)_n\subset B(x,r/2)$ such that $x_n\to x$. Since $v$ is locally bounded in $D$, there exists a constant $c_1>0$ such that $v(y)\le c_1$ for all $y\in B(x,r)$. Therefore, 
\begin{align*}
	\int_D |G_D(x_n,y)-G_D(x,y)|\rho(y)v(y)dy&\leq c_1\int_{D} |G_D(x_n,y)-G_D(x,y)|\rho(y)dy\\
	&\quad + \int_{D\cap B(x,r)^c} |G_D(x_n,y)-G_D(x,y)|\rho(y)v(y)dy
\end{align*}
Since $ G_D(x_n,y)\rho(y)\to  G_D(x,y)\rho(y)$ as $n\to\infty$, for a.e. $y\in D$, by Vitali's convergence theorem, \cite[Theorem 16.6 $(i)\iff (iii)$]{Sch}, it is enough to show that
\begin{align*}
	&\lim_{n\to \infty}\int_{D} G_D(x_n,y)\rho(y)dy=\int_{D} G_D(x,y)\rho(y)dy\quad\text{and}\\
		&\lim_{n\to \infty}\int_{D\cap B(x,r)^c} G_D(x_n,y)\rho(y)v(y)dy=\int_{D\cap B(x,r)^c} G_D(x,y)\rho(y)v(y)dy.
\end{align*}
The first limit follows directly from the assumption $G_D\rho\in C(D)$. For the second integral, we will show that there exists a constant $c_2>0$ such that
\begin{equation}\label{e:ln-1}
G_D(w,y)\leq c_2 G_D \1 (y)\, ,\quad w\in B(x,r/2), \ \ y\in D\cap B(x,r)^c.
\end{equation}
Therefore, since $\rho vG_D\1\in L^1(D)$ and $x_n\in B(x,r/2)$, we can apply the dominated convergence theorem to obtain
\begin{align*}
	\lim_{n\to\infty}\int_{D\cap B(x,r)^c} G_D(x_n,y)\rho(y)v(y)\, dy&=\int_{D\cap B(x,r)^c} G_D(x,y)\rho(y) v(y)\, dy.
\end{align*}

It remains to show \eqref{e:ln-1}. First note that $G_D(\cdot, y)$ are harmonic functions in $B(x,r)$ for all $y\in D\cap B(x,r)^c$. By the Harnack principle, there exists $c_3>0$ such that
\begin{equation}\label{e:ln-2}
G_D(w,y)\le c_3 G_D(x,y)\, ,\quad \textrm{for all }w\in B(x,r/2) \textrm{ and all }y\in D\cap B(x,r)^c. 
\end{equation}
Let $\psi:D\to [0,1]$ be a function with support in $B(x,r/2)$. Then both $G_D\psi$ and $G_D(x,\cdot)$ are regular harmonic in $D\cap B(x,r)^c$ and vanish in the sense of the limit on $D^{\mathrm{reg}}$ and by definition on $\overline{D}^c$. 

 Let $z\in \partial D$. By \cite[Theorem 1.1]{KSV16}, there exists a finite limit
$$
a(z):=\lim_{y\to z, y\in D}\frac{G_D(x,y)}{G_D\psi(y)}.
$$
Therefore, there exists a $0<\epsilon(z)<\mathrm{dist}(B(x,r), \partial D)/2$ such that
$$
\frac{G_D(x,y)}{G_D\psi(y)}\le a(z)+1\, , \quad \textrm{for all }y\in D\cap B(z,\epsilon(z)).
$$
By compactness of $\partial D$, there are finitely many points $z_1, z_2, \dots, z_n\in \partial D$ and $\delta>0$ such that $\partial D\subset \overline{D}\setminus D_{\delta}\subset \cup_{j=1}^n B(z_j, \epsilon(z_j))$. Thus for any $y\in D\setminus D_{\delta}$ it holds that
\begin{equation}\label{e:ln-3}
\frac{G_D(x,y)}{G_D\psi(y)}\le \max_{j=1, \dots, n}(a(z_j)+1)=:c_4. 
\end{equation}
Further, since both $G_D \psi$ and $G_D(x, \cdot)$ are continuous (and strictly positive) on the compact set $\overline{D_{\delta}}\cap B(x,r)^c$, we get that
\begin{equation}\label{e:ln-4}
\frac{G_D(x,y)}{G_D\psi(y)}\le c_5\, ,\quad y\in \overline{D_{\delta}}\cap B(x,r)^c.
\end{equation}
Combining \eqref{e:ln-2}--\eqref{e:ln-4} together with $G_D \psi\le G_D \1$, we get \eqref{e:ln-1}. \qed

\begin{lemma}\label{l:continuity of C0 bounded functions}
	Let $|g|\le f$ such that $G_Df\in C_0(D)$. Then $G_Dg\in C_0(D)$.
\end{lemma}
\proof
	Let $(x_n)_n\subset D$ be a sequence that converges to $x\in D$. We have
	\begin{align}\label{e:C_0 bound}
		|G_Dg(x_n)-G_Dg(x)|&\le \int_D|G_D(x_n,y)-G_D(x,y)||g(y)|dy\nonumber\\
		&\le\int_D|G_D(x_n,y)-G_D(x,y)|f(y)dy.
	\end{align}
	Since $G_D(x_n,y)f(y)\to G_D(x,y)f(y)$ as $n\to\infty$ and $G_Df\in C_0(D)$ by Vitali's theorem \cite[Theorem 16.6 $(i)\iff(iii)$]{Sch} we have that the right-hand side of \eqref{e:C_0 bound} tends to 0. Hence, $G_Dg\in C(D)$.
	
	To see that $G_Dg\in C_0(D)$ it is enough to notice that $0\le |G_Dg(x)|\le G_Df(x)$ in $D$ so when $x\to z\in\partial D$ we have $G_Dg(x)\to 0$.
	\qed


\section{The semilinear problem in bounded open set}\label{s:semilinear}

Let us now turn to the semilinear problem.
For functions $f:D\times\R\to\R$ and $u:D\to \R$ let $f_u:D\to \R$ be a function defined by
$$
f_u(x)=f(x,u(x)).
$$

\begin{defn}\label{d:wds}
Let $f:D\times\R\to\R$ be a function, $\lambda\in\mathcal M(D^c)$ and $\mu\in\mathcal M(\partial D)$
a measure concentrated on $\partial_M D$,
such that $P_D |\lambda|+M_D |\mu|<\infty$ on $D$. 
A function $u\in L^1(D)$ is called a weak dual solution to the semilinear problem 
\begin{equation}\label{e:4}
	\begin{array}{rcll}
		-Lu(x)&=& f(x,u(x))& \quad \text{in } D\\
		u&=&\lambda& \quad \text{in }D^c\\
		W_Du&=&\mu&\quad \text{on }\partial D
	\end{array}
\end{equation}
if $u$ satisfies the equality
\begin{align}
\int_D u(x)\psi(x)dx&=\int_D f(x,u(x)) G_D\psi(x) dx \nonumber \\
&+ \int_{D^c}\int_{D}P_D(x,z)\psi(x)\, dx\, \lambda(dz) \nonumber\\
& + \int_{\partial_M D}\int_{D}M_D(x,z)\psi(x)\, dx\, \mu(dz), \label{e:weak-dual}
\end{align}
for every $\psi\in C_c^\infty(D)$. If in the equation above we have $\geq$ ($\leq$) instead of the equality and the inequality holds for every nonnegative $\psi\in C_c^\infty(D)$, we say that $u$ is a supersolution (subsolution) to \eqref{e:4}.
\end{defn}
\begin{remark}\label{r:repres. of solution}
{\rm
\begin{enumerate}[leftmargin=0cm,itemindent=.5cm,labelwidth=\itemindent,labelsep=0.03cm,align=left, label=(\roman*)]
\item
Recall from Subsections \ref{ss:killed} and \ref{ss:poisson-martin} that if $P_D|\lambda|(x)+M_D|\mu|(x)<\infty$ for some $x\in D$, then $P_D|\lambda|(x)+M_D|\mu|(x)<\infty$ for all $x\in D$. 
Also, since $P_D|\lambda|<\infty$, $\lambda$ is a measure on $\R^d\setminus (D\cup \partial_MD)$, see Subsection \ref{ss:poisson-martin}, so conditions in \eqref{e:4} in $D^c$ and on $\partial D$ are indeed complementary.
\item
Note that by Fubini's theorem and symmetry of $G_D$, the above definition implies that the weak dual solution $u$ of \eqref{e:4} satisfies
\[
u(x)=G_Df_u(x)+P_D \lambda(x)+ M_D \mu(x),
\]
for almost every $x\in D$. Moreover, if we set $g=P_D\lambda + M_D\mu$, then (3.2) is equivalent to
\begin{equation}\label{e:weak-dual-alt}
\int_D u(x)\psi(x)dx=\int_Df(x,u(x))G_D\psi(x)dx + \int_D g(x)\psi(x) dx.
\end{equation}
Also, suppose that $u\in L^1_{loc}(D)$ satisfies \eqref{e:4}. This also implies that $u=G_Df_u+P_D \lambda+ M_D \mu$ a.e. in $D$. Since $G_Df_u$, $P_D \lambda$, $M_D \mu\in L^1(D)$, see Subsections \ref{ss:killed} and \ref{ss:poisson-martin}, we have $u\in L^1(D)$, i.e. every function that satisfies \eqref{e:4} must be in $L^1(D)$. 
\end{enumerate}
}
\end{remark}

Before we show an existence and uniqueness theorem for a wide class of problems we show an auxiliary result. For a Borel set $A\subset D$ and $x\in A$, let $\omega_A^x(dz):=\P_x(X_{\tau_A}\in dz)$ denote the harmonic measure. If $u:\R^d\to [-\infty, \infty]$, let $P_A u(x):=\E_x[u(X_{\tau_A})]=\int_{\R^d}u(y)\omega_A^x(dy)$ whenever the integral makes sense. We also recall that $G_A(x,y)=0$ if $y\notin A$. Finally, if the function $u$ is defined only on $D$, we extend it to all of $\R^d$ by setting $u(x)=0$ for $x\notin D$, and denote the extended function by $u \1_D$. 

\begin{lemma}\label{l:u-in-A}
	Let $D$ be an open bounded set in $\R^d$, $f:D\to[-\infty,\infty]$ a function on $D$ and  $\lambda\in\mathcal M(D^c)$ such that 
	\[
	G_D|f|(x_0), P_D|\lambda|(x_0)<\infty \text{ for some }x_0\in D. 
	\]
	Let $u$ be a function on $D$ satisfying
	\[
	u(x)=G_Df (x)+P_D\lambda(x) \text{ for a.e. }x\in D
	\] 
	and $A\subset D$ an open set. Then for a.e. $x\in A$, 
	\begin{equation}\label{eq:onA}
	u(x)= G_Af(x)+ P_A(u\1_{D})(x) +\int_{D^c} P_A(x,y)\lambda(dy).
	\end{equation}
\end{lemma}
\proof
First recall that if $G_D|f|(x_0), P_D|\lambda|(x_0)<\infty$ for some $x_0\in D$ then $G_D|f|(x), P_D|\lambda|(x)<\infty$ for almost every $x\in D$, see Subsections \ref{ss:killed} and \ref{ss:poisson-martin}. By the strong Markov property we have that
$$
G_D(x,y)=G_A(x,y)+\int_{D\setminus A}G_D(z,y)\omega_A^x(dz),\quad x\in A, y\in D,
$$
and then \eqref{e:poisson-kernel} implies that
$$
P_D(x,y)=P_A(x,y)+\int_{D\setminus A}P_D(z,y)\omega_A^x(dz),\quad x\in A,\, y\in D^c.
$$
Therefore, for a.e. $x\in A$ we have
\begin{align*}
u(x)&=\int_A G_D(x,y)f(y) dy+\int_{D\setminus A} G_D(x,y)f(y) dy+ \int_{D^c} P_D(x,y)\lambda(dy)\\
&=\int_A G_A(x,y)f(y)dy+\int_A \int_{D\setminus A}G_D(z,y)\omega_A^x(dz) f(y)dy\\
&\ \ +\int_{D\setminus A} \int_{D\setminus A}G_D(z,y)\omega_A^x(dz) f(y)dy+ \int_{D^c} P_D(x,y)\lambda(dy)\\
&=\int_A G_A(x,y)f(y) dy+\int_{D\setminus A}\left(\int_D G_D(z,y) f(y)dy\right) \omega_A^x(dz)\\
&\ \ + \int_{D^c} P_D(x,y)\lambda(dy)\\
&=\int_A G_A(x,y)f(y)dy+\int_{D\setminus A}u(z)\omega_A^x(dz)\\
&\ \ -\int_{D\setminus A}\left(\int_{D^c}P_D(z,y)\lambda(dy)\right) \omega_A^x(dz)+ \int_{D^c} P_D(x,y)\lambda(dy)\\
&=\int_A G_A(x,y)f(y)dy+\int_{D\setminus A}u(z)\omega_A^x(dz)+ \int_{D^c} P_A(x,y)\lambda(dy).
\end{align*}
\qed   
\begin{remark}\label{r:u-in-A}
	{\rm
		Let $u=G_D f+P_D\lambda$ as above and set $u=\lambda$ on $D^c$. For an open set $A\subset D$ with a Lipschitz boundary consider the linear problem $-L u_A=f$ in $A$, $u_A=u$ in $A^c$, and $W_A u_A=0$ on $\partial A$. Then Lemma \ref{l:u-in-A} says that $u_A=u$ in $A$. 
	}
\end{remark}

\begin{prop}\label{p:uniqueness}
	Let $D\subset\R^d$ be a bounded open set and let $f:D\times\R\to\R$ be a function which is nonincreasing in the second variable. Then the continuous weak dual solution to \eqref{e:4} is unique.
\end{prop}
\proof
Let $u_1$ and $u_2$ be two continuous solutions to \eqref{e:4}. Remark \ref{r:repres. of solution}(ii) yields that $u_i=G_Df_{u_i}+P_D\lambda+M_D\mu$ a.e. on $D$, $i=1,2$, hence $u_1-u_2=G_Df_{u_1}-G_Df_{u_2}$ a.e. on $D$. Note that $A:=\{x\in D: u_1(x)>u_2(x)\}$ is open and that $f(x,u_1(x))\le f(x,u_2(x))$, $x\in A$, since $f$ is nonincreasing. Using Lemma \ref{l:u-in-A} we get for a.e. $x\in A$
\begin{align*}
	0<u_1(x)-u_2(x)=G_A(f_{u_1}-f_{u_2})(x)+P_A\big((u_1-u_2)\mathbf{1}_D\big)(x)\le 0
\end{align*}
hence $A=\emptyset$. Similarly we get $\{x\in D: u_2(x)>u_1(x)\}=\emptyset$. 
\qed

Let us recall the condition \textbf{(F)} on the function $f$:

\medskip
\noindent 
\textbf{(F)} $f:D\times\R\to\R$ is continuous in the second variable and there exist a function $\rho:D\to[0,\infty)$ and a continuous function $\Lambda:[0,\infty)\to [0,\infty)$ such that $|f(x,t)|\leq \rho(x)\Lambda(|t|)$.

\begin{thm}\label{t:any-f}
Let $D\subset\R^d$ be a  
bounded open set and let $f:D\times\R\to\R$ be a function satisfying the condition \textbf{(F)}. Let $\lambda\in\mathcal M(D^c)$ such that $P_D|\lambda|<\infty$ and $\mu\in\mathcal M(\partial D)$ be a finite measure concentrated on $\partial_M D$. Assume that the nonlinear problem \eqref{e:4} admits a weak dual subsolution $\underline{u}\in L^1(D)\cap C(D)$ and a weak dual supersolution $\overline{u}\in L^1(D)\cap C(D)$ such that $\underline{u}\leq \overline{u}$. Set $g:=P_D\lambda+M_D\mu$ and $h:=|\overline{u}|\vee |\underline{u}|$. 
If one of the following conditions holds
\begin{enumerate}[leftmargin=0cm,itemindent=.5cm,labelwidth=\itemindent,labelsep=0.03cm,align=left, label=(\roman*)]
\item  $\mu\equiv 0$, $G_D\rho\in C_0(D)$ and $\underline{u}, \overline{u}\in L^\infty(D)$ such that for every open subset $A\subset D$ and a.e. $x\in A$ 
\begin{align}
	& \underline{u}(x)\leq G_A f_{\underline{u}}(x)+P_A(\underline{u}\1_{D})(x)+P_A\lambda(x),\label{eq:sub}\\ 
	&\overline{u}(x)\geq G_Af_{\overline{u}}(x)+ P_A(\overline{u}\1_{D})(x)+ P_A\lambda(x);\label{eq:super}
\end{align}
\item $\mu\equiv 0$, $\Lambda$ is nondecreasing, $G_D(\rho\Lambda(h))\in C_0(D)$ and $\underline{u}$ and $\overline{u}$ satisfy \eqref{eq:sub} and \eqref{eq:super}, respectively;
\item $\Lambda$ is nondecreasing, $G_D(\rho\Lambda(h))\in C_0(D)$ and there exists a constant $C>0$ such that, on $D$, $G_D(\rho\Lambda(h))\leq C$ and $\underline{u}-g\leq - C<C\leq \overline{u}-g$;
\end{enumerate}
then \eqref{e:4} has a weak dual solution $u\in L^1(D)\cap C(D)$ satisfying
\begin{align}\label{e:between solution}
	\underline{u}\leq u\leq \overline{u}. 
\end{align}
If, in addition, $f$ is nonincreasing in the second variable, then $u$ is a unique continuous weak dual solution to \eqref{e:4}.
\end{thm}

\begin{remark}
{\rm
Note that by Lemma \ref{l:u-in-A} a supersolution $\overline{u}$ to the nonlinear problem \eqref{e:4} satisfies the condition \eqref{eq:super} if, for example, $\overline{u}$ is a solution to the nonlinear problem 
\[	\begin{array}{rcll}
		-Lu(x)&=& f(x,u(x))& \quad \text{in } D\\
		u&=&\widetilde{\lambda}& \quad \text{in }D^c\\
		W_Du&=&0&\quad \text{on }\partial D
	\end{array}
	\]
}
for some $\widetilde\lambda\in\mathcal M(D^c)$ such that $P_D |\widetilde\lambda|<\infty$ on $D$ and $\lambda\le \widetilde\lambda$ (for details see the proof of Theorem \ref{t:nonpositive-f} and the functions $u_{n,k}$).
\end{remark}

\pff of Theorem \ref{t:any-f}.
First note that by using \eqref{e:weak-dual-alt} and \eqref{e:W_D(G_D)}, a function $u\in L^1(D)$ is the solution to \eqref{e:4} if and only if $u-g$ is the solution to the homogeneous problem
\begin{equation}\label{e:4_hom}
	\begin{array}{rcll}
		-Lw(x)&=& f(x,w+g)& \quad \text{in } D\\
		w&=&0& \quad \text{in } D^c\\
		W_Dw&=&0&\quad \text{on }\partial D.
	\end{array}
\end{equation}
Thus, we solve \eqref{e:4_hom}. For general $v\in C_0(D)$, the function $f_v$ need not satisfy the Kato condition $G_D|f_v|\in C_0(D)$,
so we define a modification of $f$ in the following way:
\begin{equation}\label{e:F}
	F(x,t)=\begin{cases}
		f(x,\overline{u}(x)), & t>\overline{u}(x)-g(x)\\
		f(x,t+g(x)), & \underline{u}(x)-g(x)\leq t\leq \overline{u}(x)-g(x)\\
		f(x,\underline{u}(x)), & t<\underline{u}(x)-g(x).
	\end{cases}
\end{equation}
Note that $F$ is continuous in the second variable. Furthermore, 
\begin{equation}\label{eq:GDF}
\textrm{if } v\in C_0(D), \textrm{ then } G_D|F_v|\in C_0(D),
\end{equation}
since 
\begin{itemize}
	\item under $(i)$, $G_D\rho\in C_0(D)$ and
	\begin{equation}\label{eq:Fu1}
	|F(x,v(x))|\leq \rho(x)\max_{y\in[0,M]}\Lambda(y),
	\end{equation}
	where $M:=\max\{\|\underline{u}\|_\infty, \|\overline{u}\|_\infty\}$ and $c_1:=\max_{y\in[0,M]}\Lambda(y)<\infty$ so the claim now follows from Lemma \ref{l:continuity of C0 bounded functions};
	\item under $(ii)$ and $(iii)$, $G_D(\rho\Lambda(h))\in C_0(D)$ and
		\begin{equation}\label{eq:Fu2}
	|F(x,v(x))|\leq \rho(x)\Lambda(|\underline{u}(x)|\vee |\overline{u}(x)|)=\rho(x)\Lambda(h(x)),
	\end{equation}
and the claim again follows from Lemma \ref{l:continuity of C0 bounded functions}.
\end{itemize}
Next we consider an auxiliary problem
\begin{equation}\label{e:5}
	\begin{array}{rcll}
		- Lu(x)&=& F(x,u)& \quad \text{in } D\\
		u&=&0& \quad \text{in } D^c\\
		W_Du&=&0&\quad \text{on }\partial D,
	\end{array}
\end{equation} 
whose solution will be given by the Schauder fixed point theorem. To this end, 
\begin{itemize}
\item under $(i)$, set $C:=\|G_D\rho\|_{L^\infty(D)}\|\Lambda\|_{L^\infty([0,M])}$;
\item under $(ii)$, set $C:=\|G_D(\rho\Lambda(h))\|_{L^\infty(D)}$;
\item under $(iii)$, let $C$ be the constant from the assumption $(iii)$;
\end{itemize} 
and let $K=\{v\in C_0(D): \|v\|_\infty\leq C\}$. 
Define the operator $T$ by
\begin{align}\label{e:operator T}
Tv(x)=\int_D F(y,v(y))G_D(x,y)dy,\ v\in C_0(D).
\end{align}
From \eqref{eq:GDF} we have $Tv\in C_0(D)$. We now prove the  continuity of $T$. Suppose the opposite, i.e.~suppose that there are $\varepsilon>0$, $(x_n)_n\subset D$, $(v_n)_n\subset C_0(D)$ and $v\in C_0(D)$ such that $||v_n-v||_\infty\to0$ and $|Tv_n(x_n)-Tv(x_n)|\ge \varepsilon$, for all $n\in\N$. Since $\overline D$ is compact there is $x\in\overline D$ and a subsequence of $(x_n)_n$ denoted again by $(x_n)_n$ such that $x_n\to x$. We have
\begin{align}\label{e:T contradict}
	\varepsilon\le |Tv_n(x_n)-Tv(x_n)|\le |Tv_n(x)-Tv(x)|+|Tv_n(x_n)-Tv_n(x)|+|Tv(x)-Tv(x_n)|.
\end{align}
Note that if $x\in \partial D$, then $Tv_n(x)=Tv(x)=0$ by \eqref{e:operator T}. Since $F$ is continuous in the second variable using the dominated convergence theorem with bounds from \eqref{eq:Fu1} and \eqref{eq:Fu2} for the first term, for $x\in D$ we have $|Tv_n(x)-Tv(x)|\to0$ as $n\to \infty$.  For the second and the third term let us also look first at the case $x\in\partial D$. Note that from \eqref{eq:Fu1} and \eqref{eq:Fu2} we have
\begin{itemize}
	\item under $(i)$
	\begin{align*}
	|Tw(x_n)|\le c_1\int_DG_D(x_n,y)\rho(y)dy=c_1G_D\rho(x_n)\to0,\enskip\text{ as $x_n\to x$,} \quad w\in\{v,v_n\},
	\end{align*}
	since $G_D\rho\in C_0(D)$;
	\item under $(ii)$ and $(iii)$
	\begin{align*}
	|Tw(x_n)|\le \int_DG_D(x_n,y)\rho(y)\Lambda(h(y))dy=G_D(\rho\Lambda(h))(x_n)\to0,\enskip\text{ as $x_n\to x$,} \quad w\in\{v,v_n\},
	\end{align*}
	since $G_D(\rho\Lambda(h))\in C_0(D)$.
\end{itemize}
If $x\in D$ then $G_D(x_n,y)\to G_D(x,y)$ so using \cite[Theorem 16.6 $(i)\iff (iii)$]{Sch}
\begin{itemize}
	\item under $(i)$
	\begin{align*}
	|Tw(x_n)-Tw(x)|\le c_1\int_D|G_D(x_n,y)-G_D(x,y)|\rho(y)dy\to0,\enskip\text{ as $x_n\to x$,} \quad w\in\{v,v_n\},
	\end{align*}
	since $G_D\rho\in C_0(D)$;
	\item under $(ii)$ and $(iii)$
	\begin{align*}
	|Tw(x_n)-Tw(x)|\le \int_D|G_D(x_n,y)-G_D(x,y)|\rho(y)\Lambda(h(y))dy\to0,\enskip\text{ as $x_n\to x$,} \quad w\in\{v,v_n\},
	\end{align*}
	since $G_D(\rho\Lambda(h))\in C_0(D)$.
\end{itemize}
Thus, we have a contradiction with \eqref{e:T contradict}, i.e. $T$ is continuous.

Also, from \eqref{eq:Fu1}, \eqref{eq:Fu2} and the choice of constant $C$ we get $T(K)\subset K$.

We are left to prove that $T(K)$ is a precompact subset of $K$. By Arzel\`{a}-Ascoli theorem it suffices to note that the functions $\{Tv : v \in K\}$ are equicontinuous by the same calculations as above.

Hence by the Schauder fixed point theorem there is a function $u\in K$ such that  
\[
	u(x)=\int_D F(y,u(y))G_D(x,y)dy,
\]
i.e. $u$ is a weak dual solution to \eqref{e:5}. It follows immediately from \eqref{e:F} that, if $\underline{u}-g\leq u\leq \overline{u}-g$, then $u$ is also a  weak dual solution to \eqref{e:4_hom}. Finally, we show that the obtained solution $u$ to \eqref{e:5} is between $\underline{u}-g$ and $\overline{u}-g$. In case of assumption (iii), this is obvious. Under (i) or (ii), set $A=\{x\in D:u(x)>\overline{u}(x)-g(x)\}$. Note that $F_u(y)=f_{\overline{u}}(y)$ for all $y\in A$ and that $A$ is an open subset of $D$, since both $u$ and $\overline{u}-g$ are continuous on $D$. Then, for every $x\in A$, by \eqref{eq:onA} we have 
\begin{align*}
	u(x)+g(x)&=G_AF_u(x) + P_A ((u+g)\1_D)(x)+P_A\lambda(x)\\
	&\leq G_Af_{\overline u}(x)+ P_A(\overline{u}\1_D)(x)+P_A\lambda(x)\\
	&\leq\overline{u}(x),
\end{align*}
where the first inequality comes only from the middle term and the second one is \eqref{eq:super}. This implies that $A=\emptyset$. By using \eqref{eq:sub}, one can analogously show that $\{x\in D:u(x)\leq\underline{u}(x)-g(x)\}=\emptyset$.

Uniqueness follows from Proposition \ref{p:uniqueness}.
\qed

In the following corollary we extend the main result from \cite{BJK} to our setting of more general non-local operators.

\begin{corollary}\label{c:bogdan}
Let $D\subset\R^d$ be a  bounded open set and let $f:D\times\R\to\R$ be a function satisfying the condition \textbf{(F)} with $\Lambda$ nondecreasing. Let $\lambda\in\mathcal M(D^c)$  such that $P_D|\lambda|<\infty$ and $\mu\in\mathcal M(\partial D)$ a finite measure concentrated on $\partial_M D$. Set $g:=P_D\lambda+M_D\mu$ and $\overline{g}:=P_D|\lambda| + M_D |\mu|$. Assume that $G_D\rho \in C_0(D)$, $G_D(\rho\Lambda(2\overline{g}))\in C_0(D)$, and that either (a) $\Lambda$ is sublinearly increasing, $\lim_{t\to \infty}\Lambda(t)/t=0$, or (b) $m$ is sufficiently small. Then the semilinear problem
\begin{equation}\label{e:bogdan-1}
\begin{array}{rcll}
		- Lu(x)&=& m f(x,u(x))& \quad \text{in } D\\
		u&=&\lambda& \quad \text{in } D^c\\
		W_Du&=&\mu&\quad \text{on }\partial D
\end{array}
\end{equation}
has a weak dual solution $u\in L^1(D)\cap C(D)$ such that $|u|\leq \overline{g}+C$, for some $C>0$.

If, in addition, $f$ is nonincreasing in the second variable, $u$ is a unique continuous weak dual solution to \eqref{e:bogdan-1}.
\end{corollary}
\proof
We use Theorem \ref{t:any-f}(iii)  with $m f$ instead of $f$ and first choose the constant $C>0$. Set $r_1:=\sup_{x\in D} G_D\rho(x)$ and $r_2:=\sup_{x\in D} G_D (\rho\Lambda(2\overline{g}))(x)$. By the assumption, we have that $r_1<\infty$ and $r_2<\infty$. If (b) holds, given any $C>0$ we can find $m$ small enough such that $m(\Lambda(2C)r_1+r_2)\le C$. If (a) holds, then since $\Lambda$ is sublinearly increasing, we can find $C>0$ large enough so that again $m(\Lambda(2C)r_1+r_2)\le C$. 

Let $\overline{u}:=C+\overline{g}$, $\underline{u}:=-\overline{u}$ and $h:=|\overline{u}|\vee |\underline{u}|=C+\overline g$. Clearly, $\overline{u}$ and $\underline{u}$ belong to $L^1(D)\cap C(D)$ and satisfy $\underline{u}-g\leq - C<C\leq \overline{u}-g$. We check that $\overline{u}$ is a supersolution  of \eqref{e:bogdan-1}. Indeed, 
\begin{eqnarray}
\lefteqn{|G_D(m f_{\overline{u}})+g|\le mG_D|f_{C+\overline{g}}|+\overline{g}\le mG_D(\rho\Lambda(C+\overline g))+\overline{g}}\nonumber\\
&\le &m G_D\big(\rho \big(\Lambda(2C)+\Lambda(2\overline{g})\big)\big) +\overline{g}\le m\big(\Lambda(2C)r_1+r_2\big)+\overline{g}\le C+\overline{g}=\overline{u}.\nonumber
\end{eqnarray}
In the same way we see that $\underline{u}$ is a subsolution. It remains to check that  $G_D(m\rho \Lambda(h))\in C_0(D)$ and $G_D(m\rho \Lambda(h))\le C$. By the same computations as above we have 
\begin{align}
G_D(m\rho \Lambda(h))&\le m\Lambda(2C)G_D\rho+mG_D\rho\Lambda(2\overline{g}))\label{e: Kato}\\
&\leq m(\Lambda(2C)r_1+r_2)\le C.\nonumber
\end{align}
Since $G_D\rho\in C_0(D)$ and $G_D(\rho \Lambda(2\overline{g}))\in C_0(D)$, by \eqref{e: Kato} and Lemma \ref{l:continuity of C0 bounded functions} we also have  $G_D(m\rho \Lambda(h))\in C_0(D)$.

Uniqueness follows from Proposition \ref{p:uniqueness}.
\qed

Our next goal is to extend Corollary \ref{c:bogdan}  to a wider class of nonpositive functions $f$. First we show an additional auxiliary result. This result provides an approximation of a  nonnegative harmonic function on $D$ by an increasing sequence of potentials. It is a consequence of a rather well-known fact that we prove in the appendix, see Proposition \ref{p:approximation}. We can use this result because the semigroup $(P_t^D)_{t\ge 0}$ is strongly Feller, the process $X^D$ is transient, nonnegative harmonic functions are excessive, and the potential $G_D \1$ is continuous and satisfies $0<G_D\1 <\infty$ on $D$.
\begin{lemma}\label{l:approximation}
Let $h:D\to [0,\infty)$ be a harmonic function with respect to the process $X^D$. There exists a sequence $(\wt{f}_k)_{k\ge 1}$ of nonnegative, bounded and continuous functions such that $G_D\wt{f}_k \uparrow h$. 
\end{lemma}

\begin{thm}\label{t:nonpositive-f}
Let $D\subset\R^d$ be a bounded open set. Let $f:D\times\R\to(-\infty,0]$ be a function that satisfies \textbf{(F)} with $G_D\rho\in C_0(D)$. Assume, additionally, that $f(x,0)=0$.
Let $\lambda\in\mathcal M(\R^d\setminus \overline{D})$ be a nonnegative measure such that $P_D\lambda<\infty$ and $\mu\in \mathcal M(\partial D)$ be a finite nonnegative measure concentrated on $\partial_M D$. Let $g:=P_D\lambda+M_D\mu$. If the semilinear  problem \eqref{e:4} satisfies one of the following conditions:
	\begin{enumerate}[leftmargin=0cm,itemindent=.5cm,labelwidth=\itemindent,labelsep=0.03cm,align=left, label=(\roman*)]
		\item
		$\mu\equiv 0$;
		\item
		$\mu\not\equiv0$, the function $\Lambda$ is nondecreasing and $\rho\Lambda(g)G_D\1\in L^1(D)$;
	\end{enumerate}
 then the problem \eqref{e:4} has a  nonnegative  weak dual solution $u\in L^1(D) \cap C(D)$.  If, in addition, $f$ is nonincreasing in the second variable, then $u$ is a unique continuous  solution to \eqref{e:4}.
\end{thm}
\proof
Let $(\wt{f}_k)_k$ be a sequence of nonnegative, bounded and continuous functions on $D$ from Lemma \ref{l:approximation} such that $G_D\wt{f}_k\uparrow M_D\mu$. Let $(K_n)_n$ be an increasing sequence of compact sets such that $K_n\uparrow \overline{D}^c$. Then, for $n\in\mathbb N$ the measure $\lambda_n(\cdot)=\lambda(\cdot\cap K_n)$ is a finite nonnegative measure on $\overline{D}^c$. Consider the following semilinear problem
\begin{equation}\label{e:8}
\begin{array}{rcll}
- Lu(x)&=& f(x,u(x))+\wt{f}_k(x)& \quad \text{in } D\\
u&=&\lambda_n& \quad \text{in }D^c\\
W_Du&=&0&\quad \text{on }\partial D.
\end{array}
\end{equation}
Since $f(x,0)=0$ and $\wt f_k\geq 0$, $\underline{u}\equiv0$ is a subsolution to \eqref{e:8}. Furthermore, since $f$ is nonpositive, as a supersolution to \eqref{e:8} we take the solution $u_k^{(n)}=G_D\wt{f}_k+P_D\lambda_n$ of the linear problem
\begin{equation*}
\begin{array}{rcll}
-  Lu(x)&=&\wt{f}_k(x)& \quad \text{in } D\\
u&=&\lambda_n& \quad \text{in }D^c\\
W_Du&=&0&\quad \text{on }\partial D.
\end{array}
\end{equation*}
Fix $k\in\mathbb N$. Notice that $u^{(n)}_k\in C(D)$ and that, by Lemma \ref{l:u-in-A}, $u^{(n)}_k$ satisfies \eqref{eq:onA}. Moreover, since $\lambda_n$ is finite and
\[
\sup_{x\in D, z\in K_n}P_D(x,z)\leq j(\text{dist}(D,K_n))\sup_{x\in D}G_D\1(x)<\infty,
\]
$u_k^{(n)}$ is bounded. This means that we can apply Theorem \ref{t:any-f}(i) so that for $n=1$ the problem \eqref{e:8} has a solution $u_{1,k}\in C(D)\cap L^\infty(D)$ such that $0\leq u_{1,k}\leq u^{(1)}_k$. Note that since $\lambda_1\le \lambda_2$, $u_{1,k}$ is also a subsolution to the problem \eqref{e:8} for $n=2$ such that \eqref{eq:sub} holds for every open subset $A\subset D$, that is for a.e. $x\in A$
\begin{align*}
u_{1,k}(x)&= G_Af_{u_{1,k}}(x)+G_A\widetilde f_k(x)+P_A(u_{1,k}\1_{D})(x)+P_A\lambda_1(x)\\
	&\leq  G_Af_{u_{1,k}}(x)+G_A\widetilde f_k(x)+P_A(u_{1,k}\1_{D})(x)+P_A\lambda_2(x).
\end{align*}
Since $u_{1,k}\leq u^{(1)}_k\leq u^{(2)}_k$, again by Theorem \ref{t:any-f}(i), there exists a solution $u_{2,k}\in C(D)\cap L^\infty(D)$ to the problem \eqref{e:8} with $\lambda_2$ on $D^c$, such that $u_{1,k}\leq u_{2,k} \leq u^{(2)}_k$. By iterating this procedure, we obtain an increasing sequence $(u_{n,k})_{n\in\N}$ of solutions to problems \eqref{e:8} for different $n\in\mathbb N$. Moreover, the sequence $(u_{n,k})_{n\in\N}$ is dominated by the function $u^0_k$ associated with the linear problem
\begin{equation*}
\begin{array}{rcll}
-  Lu^{0}_k(x)&=&\wt{f}_k(x)& \quad \text{in } D\\
u^{0}_k&=&\lambda& \quad \text{in }D^c\\
W_Du^{0}_k&=&0&\quad \text{on }\partial D.
\end{array}
\end{equation*}
Hence, the pointwise limit $\lim_{n\to\infty} u_{n,k}=u_k$ is well defined in $D$. We will now show that $u_k$ is a weak dual solution to the problem
\begin{equation}\label{e:9}
\begin{array}{rcll}
-  Lu(x)&=&f(x,u(x))+\wt{f}_k(x)& \quad \text{in } D\\
u&=&\lambda& \quad \text{in }D^c\\
W_Du&=&0&\quad \text{on }\partial D.
\end{array}
\end{equation}
Take any $\psi\in C_c^\infty(D)$, $\psi\ge 0$. Then by Fatou's lemma and the continuity of  the  function $f$ in the second variable, we get that
\begin{align*}
	-\int_D f(x,u_k(x))G_D\psi(x)dx&\le -\limsup_{n\to\infty}\int_D f(x,u_{n,k}(x))G_D\psi(x)dx\\
	&=-\limsup_{n\to\infty}\int_D u_{n,k}(x)\psi(x)dx+\int_D\wt{f}_k(x)G_D\psi(x)dx\\
	&\quad+\int_D P_D\lambda(x)\psi(x)dx\\
&=-\int_D u_{k}(x)\psi(x)dx+\int_D\wt{f}_k(x)G_D\psi(x)dx+\int_D P_D\lambda(x)\psi(x)dx,
\end{align*}
where we used the monotone convergence theorem in the last line. The inequality above implies that $u_k$ is a weak dual subsolution to \eqref{e:9}. To show that $u_k$ is also a supersolution of the same problem, set $D'=\supp \psi\subset\subset D$ and build a sequence $(D_l)_{l\in \N}$ of sets with Lipschitz boundaries such that $D'\subset \subset D_l \subset \subset D$ and $D_l\uparrow D$. Obviously, $\psi \in C_c^\infty(D_l)$, and both $G_{D_l}\psi\uparrow G_D\psi$ and $P_{D_l}\lambda\uparrow P_{D}\lambda$ pointwise in $D$. Also, notice that $u^0_k=G_D\wt f_k+P_D\lambda$ is continuous, hence locally bounded. Furthermore, in $D_l$ we have
$$|f(x,u_{n,k}(x))|G_{D_l}\psi(x)\le C\rho(x)G_{D_l}\psi(x),$$
where $C:=\max_{y\in D_l}\Lambda(u^0_k(y))<\infty$, 
and $\rho G_{D_l}\psi\in L^1(D)$ since $\int_D \rho G_{D_l}\psi=\int_D\psi G_{D_l}\rho\le \int_D\psi G_D\rho<\infty$. By using the dominated convergence theorem in the first equality and Lemma \ref{l:u-in-A} in the second, we have
\begin{align*}
	\int\limits_{D_l}[f(x,u_k(x))&+\wt f_k(x)]G_{D_l}\psi(x)dx=\lim_{n\to\infty}\int\limits_{D_l}[f(x,u_{n,k}(x))+\wt f_k(x)]G_{D_l}\psi(x)dx\\
	&=\lim_{n\to \infty}\left(\int\limits_{D_l}u_{n,k}(x)\psi(x)dx-\int\limits_{D_l}P_{D_l}u_{n,k}(x)\psi(x)dx-\int\limits_{D_l}P_{D_l}\lambda_n(x)\psi(x)dx\right)\\
	&\le \lim_{n\to \infty}\left(\int\limits_{D_l }u_{n,k}(x)\psi(x)dx-\int\limits_{D_l}P_{D_l}\lambda_n(x)\psi(x)dx\right)\\
	&=\int\limits_{D_l}u_{k}(x)\psi(x)dx+\int\limits_{D_l}P_{D_l}\lambda(x)\psi(x)dx.
\end{align*}
Letting $l\to\infty$ we obtain
\begin{align*}
	\int\limits_D[f(x,u_k(x))&+\wt f_k(x)]G_{D}\psi(x)dx\le \int\limits_{D}u_{k}(x)\psi(x)dx+\int\limits_{D}P_{D}\lambda(x)\psi(x)dx,
\end{align*}
which proves that $u_k$ is a supersolution, and therefore the solution to \eqref{e:9}. Notice that for $\mu\equiv0$ we have $\wt f_k\equiv 0$ so we have found a solution to the problem \eqref{e:4} under  the  assumption $(i)$.

Suppose that we have a function $\Lambda$ with properties as in the assumption $(ii)$ of this theorem. With the Arzel\`{a}-Ascoli theorem we will now find a suitable subsequence of $(u_k)_k$ that converges to a function $u$ that is a solution to the problem \eqref{e:4}. To this end first notice that $u_k$ is given by
\begin{align}
	u_k(x)&=\int_D G_D(x,y)[f(y,u_k(y))+\wt f_k(y)]dy+\int_{D^c} P_D(x,y)\lambda(dy)\nonumber\\
	&=\int_D G_D(x,y)f(y,u_k(y))dy+G_D\wt f_k(x)+P_D\lambda(x).\label{e:12}
\end{align}
Since $f$ is nonpositive, $u_k\le g=P_D\lambda+M_D\mu$ so we have the pointwise boundedness of the family $(u_k)_k$. Since $G_D\wt f_k$ increases to the continuous function $M_D\mu$, by Dini's theorem the convergence is locally uniform so the usual $3\varepsilon$-argument gives equicontinuity of the family $(G_D \wt{f}_k)_k$ at every point $x\in D$. Also, $P_D\lambda$ is continuous in $D$ so it remains to analyse the first term. We have
\begin{align*}
	&\left|\int_D G_D(x,y)f(y,u_k(y))dy-\int_D G_D(z,y)f(y,u_k(y))dy \right| \\
			&\le \int_D |G_D(x,y)-G_D(z,y)|\rho(y)\Lambda(u_k(y))dy\\
				&\le \int_D |G_D(x,y)-G_D(z,y)|\rho(y)\Lambda(g(y))dy.
\end{align*}
Equicontinuity of the first term in \eqref{e:12} now follows from Lemma \ref{l: lemma_nakon}. Now by Arzel\`{a}-Ascoli theorem we extract a subsequence $(u_{k_l})_l$ which converges pointwise to a continuous function $u$. Without loss of generality, assume that $u_k \to u$. It remains to prove that $u$ is a weak solution of \eqref{e:4}, i.e., for every $\psi \in C_c^\infty(D)$
\begin{align}\label{e:10}
	\int_D u(x) \psi(x)dx=\int_D f(x,u(x))G_D\psi(x)dx+\int_D P_D\lambda(x)\psi(x)dx+\int_D M_D\mu(x)\psi(x)dx.
\end{align}
We know that $u_k$ satisfies
\begin{align}\label{e:11}
\int_D u_k(x) \psi(x)dx=\int_D f(x,u_k(x))G_D\psi(x)dx+\int_D P_D\lambda(x)\psi(x)dx+\int_D G_D\wt{f}_k (x)\psi(x)dx.
\end{align}
Since $u_k \to u$ pointwise and $u_k\leq g$, by the dominated convergence theorem the left-hand side of \eqref{e:11} converges to the left-hand side of \eqref{e:10}. Furthermore, by the monotone convergence theorem the last term of \eqref{e:11} converges to the last term of \eqref{e:10}. To show  the convergence of the first term on the right-hand side, note that
\begin{align*}	
	|f(x,u_k(x))G_D\psi(x)|\le c_1 \rho(x)\Lambda(g(x))G_D\1(x).
\end{align*}
Now  the assumption $(ii)$ implies boundedness in $L^1(D)$, so the convergence follows from the dominated convergence theorem. Hence, $u$ is a solution to the problem \eqref{e:4}.
 Uniqueness follows from Proposition \ref{p:uniqueness}.
\qed

\begin{remark}
{\rm
\begin{enumerate}[leftmargin=0cm,itemindent=.5cm,labelwidth=\itemindent,labelsep=0.03cm,align=left, label=(\roman*)]
\item Note that the condition $\rho\Lambda(g)G_D\1\in L^1(D)$ from Theorem \ref{t:nonpositive-f} is weaker than the condition $G_D(\rho\Lambda(2\overline{g}))\in C_0(D)$ from Corollary \ref{c:bogdan}. 
\item Recall that if $D$ is regular then $q\in\JJ$ if and only if $G_D|q|\in C_0(D)$. Hence, if we assume that $D$  is regular in Theorem \ref{t:any-f} then we can equivalently assume $\rho\in \JJ$ and $\rho\Lambda(h)\in\JJ$ instead of $G_D\rho\in C_0(D)$ and $G_D(\rho\Lambda(h))\in C_0(D)$, respectively. Obviously, a similar argument applies to Corollary \ref{c:bogdan} and Theorem \ref{t:nonpositive-f}.
\end{enumerate}	
}
\end{remark}

\section{Auxiliary results in bounded $C^{1,1}$ open sets}\label{s:aux-C11}

\subsection{The renewal function} \label{ss:renewal}
We start this section by introducing a function which plays a prominent role in studying the boundary behavior in $C^{1,1}$ open sets. 

Let $Z=(Z_t)_{t\ge 0}$ be a one-dimensional subordinate Brownian motion with the characteristic exponent $\phi(\theta^2)$, $\theta\in \R$. We can think of $Z$ as one of the components of the process $X$. Let $M_t:=\sup_{0\le s\le t} Z_s$ be the supremum process of $Z$ and let $L=(L_t)_{t\ge 0}$ be the local time of $M_t-Z_t$ at zero. We refer the readers to \cite[Chapter VI]{Ber} for details. The inverse local time $L_t^{-1}:=\inf\{s>0: L_s>t\}$ is called the ascending ladder time process of $Z$. Define the ascending ladder height process $H=(H_t)_{t\ge 0}$ of $Z$ by $H_t:=M_{L_t^{-1}}=Z_{L_t^{-1}}$ if $L_t^{-1}<\infty$ and $H_t=\infty$ otherwise. The renewal function of the process $H$ is defined as
$$
V(t):=\int_0^{\infty}\P(H_s\le t)\, ds, \quad t\in \R.
$$
Then  $V(t)=0$ for $t<0$, $V(0)=0$, $V(\infty)=\infty$, and $V$ is strictly increasing. The importance of the renewal function $V$ lies in the fact that $V_{|(0,\infty)}$ is harmonic with respect to the killed process $Z^{(0,\infty)}$. This fact was for the first time used in \cite{KSV12b} in order to obtain the precise rate of decay of harmonic functions of $d$-dimensional subordinate Brownian motion. 

In the case of the isotropic $\alpha$-stable process, it holds that $V(t)=t^{\alpha/2}$. In general, the function $V$ is not known explicitly, but under the weak scaling condition \textbf{(H)} it is known, see e.g. \cite{KSV12b}, that there is a constant $C=C(R_0)\ge 1$ such that
\begin{equation}\label{e:V-phi}
C^{-1}\phi(t^{ -2})^{-1/2}\le V(t) \le C\phi(t^{ -2})^{-1/2}\, ,\quad 0<t<R_0.
\end{equation}
For more general results, covering also the case
$R_0 = \infty$, see \cite[Theorem 4.4 and Remark 4.7]{KMR3}.

Note that \eqref{e:V-phi} and weak scaling \eqref{e:wsc} of $\phi$ imply that  for all $R_1\ge 1$ there are constants  $0<\wt{a}_1\le \wt{a}_2$ depending on $R_1$ such that
\begin{equation}\label{e:wsc-V}
\wt{a}_1 \left(\frac{t}{s}\right)^{\delta_1} \le \frac{V(t)}{V(s)} \le \wt{a}_2 \left(\frac{t}{s}\right)^{\delta_2} , \quad 0<s\le t\le R_1.
\end{equation}

\subsection{Estimates in $C^{1,1}$ open set}\label{ss:estimates}
Recall that an open set $D$ in $\bR^d$ ($d\ge 2$) is said to be a $C^{1,1}$ open set if there exist a localization radius $R>0$ and a constant $\Lambda>0$ such that for every $z\in\partial D$, there exist a $C^{1,1}$ function $\psi=\psi_z: \bR^{d-1}\to \bR$ satisfying $\psi (0)=0$,  $\nabla\psi (0)=(0, \dots, 0)$, $\| \nabla \psi \|_\infty \leq \Lambda$, $| \nabla \psi (x)-\nabla \psi (z)| \leq \Lambda |x-z|$, and an orthonormal coordinate system $CS_z$: $y=(y_1, \cdots, y_{d-1}, y_d):=(\wt y, \, y_d)$ with origin at $z$
such that
$$
B(z,R)\cap D=\{ y=(\wt y, \, y_d)\in B(0, R) \mbox{ in } CS_z: y_d
>\psi (\wt y) \}.
$$
The pair $(R, \Lambda)$ is called the characteristics of the $C^{1,1}$ open set $D$. We remark that in some literature, the $C^{1,1}$ open set defined above is called a  uniform $C^{1,1}$ open set since $(R, \Lambda)$ is universal for all $z\in \partial D$.

From now until the end of this section let $D$ be a bounded open $C^{1,1}$ set. It is well known that all boundary points of a $C^{1,1}$ open set are regular and accessible. Thus, $\partial_M D=\partial D$.  Recall that  $\delta_D(x)$ denotes the distance of the point $x\in D$ to the boundary $\partial D$, while $\delta_{D^c}(z)$ denotes the distance of $z\in \overline{D}^c$ to $\partial D$. 

Under the weak scaling condition \textbf{(H)}  the following sharp two-sided estimates of the Green function, Martin kernel and the Poisson kernel are known. The comparability constant depends on the constants in \eqref{e:wsc} and the diameter of $D$. We give the estimates in terms of the renewal function $V$:
\begin{eqnarray}
G_D(x,y)&\asymp & \left( 1\wedge \frac{V(\delta_D(x))}{V(|x-y|)}\right) \left( 1\wedge \frac{V(\delta_D(y))}{V(|x-y|)}\right)\frac{V(|x-y|)^2}{|x-y|^d}\, ,\quad x,y\in D, \label{e:G_D-estimate}\\
M_D(x,z)&\asymp &\frac{V(\delta_D(x))}{|x-z|^d}\, , \quad x\in D, z\in \partial D, \label{e:M_D-estimate}\\
P_D(x,z)&\asymp & \frac{V(\delta_D(x))}{V(\delta_{D^c}(z))(1+V(\delta_{D^c}(z)))}\frac{1}{|x-z|^d}\, , \quad x\in D, z\in \overline{D}^c. \label{e:P_D-estimate}
\end{eqnarray}
For \eqref{e:G_D-estimate} see \cite[Theorem 7.3(iv)]{CKS}, \eqref{e:M_D-estimate} follows immediately from \eqref{e:martin-kernel} and \eqref{e:G_D-estimate}, while \eqref{e:P_D-estimate} is proved in \cite[Theorem 1.3]{KK}. We will also need sharp two-sided estimates of the killing function 
\begin{equation}\label{e:killing}
\kappa_D(x):=\int_{D^c}j(|y-x|)\, dy,\quad x\in D.
\end{equation}
It holds that
\begin{equation}\label{e:killing-fn-estimate}
\kappa_D(x)\asymp V(\delta_D(x))^{-2}\, ,\qquad x\in D.
\end{equation}
The upper bound is straightforward and valid in any open set $D$, while the lower bound holds in open sets satisfying the outer cone condition, see e.g.~\cite[proof of Lemma 5.7]{KSV16b}.

\subsection{Green and Poisson potentials}\label{ss:g-p-potentials}
In this subsection we state two results which should be of independent interest. The first one gives sharp two-sided estimates of the Green potential of the function $x\mapsto U(\delta_D(x))$ for a function $U:(0,\infty)\to [0,\infty)$ satisfying certain assumptions. The estimates are given in terms of the function $U$ and the renewal function $V$. A similar result was shown in \cite[Theorem 3.4]{AGCV19}. Since our proof is modeled after and is very similar to the one in \cite{AGCV19}, we defer the proof to Appendix. The second result is a sort of a counterpart of the first one and gives sharp two sided estimates of the Poisson potential of the function $z\mapsto \wt{U}(\delta_{D^c}(z))$ for a function $\wt{U}:(0,\infty)\to [0, \infty)$. The proof of this second result is simpler and will be also given in Appendix. 

To be more precise, let $U:(0,\infty)\to [0, \infty)$ be a function satisfying the following conditions:

\begin{itemize}
	\item[(U1)] Integrability condition: It holds that 
		\begin{equation}\label{e:int-cond}
		\int_{0}^1U(t)V(t)\, dt <\infty;
		\end{equation}

	\item[(U2)]  Almost nonincreasing condition: There exists  $C>0$  such that
		 \begin{equation}\label{e:wus-U}
		U(t)\le C U(s), \quad 0<s\le t\le 1;
		\end{equation}

	\item[(U3)] Reverse doubling condition: There exists $C>0$   such that
		\begin{equation}\label{e:doubling-U}
		U(t)\le C U(2t),\quad t\in (0,1); 
		\end{equation}
		\item[(U4)] Boundedness away from zero: $U$ is bounded from above on $[c, \infty)$ for each $c >0$.
\end{itemize}
We will refer to (U1)--(U4) as conditions \textbf{(U)}. Note that if $U(t)=t^{-\beta}$, $\beta\in \R$, satisfies \eqref{e:int-cond}, then it satisfies \textbf{(U)}. In particular, if the process $X$ is isotropic  $\alpha$-stable, then \eqref{e:int-cond} (hence \textbf{(U)}) is equivalent to $-\beta+\alpha/2>-1$. 

\begin{prop}\label{p:green-potential-estimate}
Assume that a function $U:(0, \infty)\to [0,\infty)$ satisfies conditions \textbf{(U)}. Then
\begin{equation}\label{e:green-potential-estimate}
G_D(U(\delta_D))(x)\asymp \frac{V(\d(x))}{\d(x)}\int_0^{\d(x)}U(t)V(t)\, dt + V(\d(x))\int_{\d(x)}^{\mathrm{diam}(D)} \frac{U(t)V(t)}{t}\, dt\, .
\end{equation}
Morover, if $U$ is positive and bounded on every bounded  subset of $(0,\infty)$, then
\begin{equation*}
G_D(U(\delta_D))(x)\asymp V(\d(x)).
\end{equation*}

\end{prop}
The asymptotic behavior of $G_D(U(\d))$ is given by the  largest term  that appears in \eqref{e:green-potential-estimate}. In this generality, this is not easy to determine (but see \cite[Theorem 3.4]{AGCV19}). It will follow from the proof that $G_D(U(\delta_D))<\infty$ if and only if \eqref{e:int-cond} holds true. Clearly, if $f:D\to [0,\infty)$ is such that $f(x)\asymp U(\d(x))$, then $G_D f(x)$ is asymptotically equal to the right-hand side of \eqref{e:green-potential-estimate}.

\begin{prop}\label{p:poisson-potential-estimate}
Let  $g:\overline{D}^c\to [0,\infty)$ be such that 
\begin{equation}\label{e:g_G}
g(y)\asymp \widetilde U(\dc(y)),\ y\in \overline{D}^c,
\end{equation}
holds for some function $\widetilde U:(0,\infty)\to [0,\infty)$. Assume that $\wt{U}$ is  bounded on every compact subset of $(0,\infty)$  and satisfies
\begin{equation}\label{e:poisson-finite}  
\int_{0}^1  \frac {\widetilde U(t)}{V(t)}dt +\int_1^\infty\frac {\widetilde U(t)}{V(t)^{2}t}dt<\infty\, .
\end{equation}
Then 
\begin{equation}\label{e:P_Dg}
P_Dg(x)\asymp V(\delta_D(x)) \int_0^{\mathrm{diam}(D)}\frac{\widetilde U(t)}{V(t)(\d(x)+t)}dt,\ x\in D,
\end{equation}
and 
\begin{equation}\label{e:P_Dg2}
P_Dg(x)\preceq\frac{V(\delta_D(x))}{\delta_D(x)},\ x\in D.
\end{equation}
\end{prop}
\begin{remark}\label{r:poisson-potential-estimate}
{\rm
In the case of the fractional Laplacian and the power function $\wt{U}(t)=t^{-\beta}$, condition \eqref{e:poisson-finite} becomes $-\alpha<\beta <1-\alpha/2$. Further, it is easy to see that for $-\beta < \alpha/2$, the integral in \eqref{e:P_Dg} is comparable to $\d(x)^{-\beta-\alpha/2}$, in the case  $\beta=-\alpha/2$ it is comparable to $\log(1/\d(x))$, while for $-\beta>\alpha/2$ it is comparable to a constant. We conclude that for  $g(y)=\delta_{D^c}(y)^{-\beta}$
$$
P_D g(x)\asymp \left\{\begin{array}{ll}
\d(x)^{-\beta}, & -\alpha<\beta<-\alpha/2, \\
\d(x)^{\alpha/2}\log(1/\d(x)),& \beta=-\alpha/2, \\
\d(x)^{\alpha/2}, & -\alpha/2<\beta<1-\alpha/2. \end{array}\right.
$$

}
\end{remark}

\subsection{Boundary estimates of harmonic functions}
Let $\sigma$ denote the $(d-1)$-dimensional Hausdorff measure on $\partial D$. It follows immediately from \eqref{e:M_D-estimate} and the estimate
$$
\int_{\partial D}\frac{\sigma(dz)}{|x-z|^d}\asymp \frac{1}{\d(x)}\, ,\quad x\in D,
$$
that
\begin{equation}\label{e:M_D-sigma}
M_D \sigma (x)=\int_{\partial D} M_D(x,z)\sigma(dz)\asymp \frac{V(\d(x))}{\d(x)}\, ,\quad x\in D.
\end{equation}

The following result appears as \cite[Theorem 4.2]{BD} for the fractional Laplacian.

\begin{prop}\label{p:boundary-limit}
Let $h\in L^1(\partial D, \sigma)$ and let $\mu(d\zeta)=h(\zeta)\sigma(d\zeta)$. If $h$ is continuous at $z\in \partial D$, then
\begin{equation}\label{e:boundary-limit-2}
\lim_{x\to z, x\in D}\frac{M_D \mu(x)}{M_D \sigma(x)}=h(z).
\end{equation}
\end{prop}
Since the proof is essentially the same as the proof of \cite[Theorem 4.2]{BD}, we omit it. Proposition \ref{p:boundary-limit}  has the following two consequences. Assume that $h$ is nonnegative, continuous in $D$, not identically equal to zero, and set $\mu(d\zeta)=h(\zeta)\sigma(d\zeta)$. Then since both $M_D\mu$ and $M_D\sigma$ are continuous and $D$ is bounded, we first conclude that there exists  $C=C(h)>0$  such that
$$
M_D\mu(x)\le C M_D\sigma(x), \quad x\in D.
$$
Secondly, there exist $z\in \partial D$, $\epsilon >0$, and  $C=C(h)>0$ such that
$$
M_D\mu(x)\ge C M_D \sigma(x), \quad x\in D\cap B(z, \epsilon).
$$
Together with \eqref{e:M_D-sigma}, these last two estimates imply that there is a constant  $C=C(h)>1$  such that
\begin{eqnarray}
M_D\mu(x)&\le & C \frac{V(\delta_D(x))}{\delta_D(x)}, \quad x\in D, \label{e:upper-bound-M}\\
M_D\mu(x)&\ge & C^{-1} \frac{V(\delta_D(x))}{\delta_D(x)}, \quad x\in D\cap B(z, \epsilon). \label{e:lower-bound-M}
\end{eqnarray}

\subsection{Kato class revisited}\label{ss:kato-revisited}
In this subsection we give a sufficient condition for a function of the distance to the boundary to be in the Kato class $\JJ$. First, note that by \eqref{e:poisson-kernel} and \eqref{e:killing}, we have that 
\begin{equation}\label{e:G-kappa-le-1}
\sup_{x\in D} G_D \kappa_D(x)\le 1.
\end{equation}
Recall from \eqref{e:killing-fn-estimate} that $\kappa_D(x)\asymp V(\d(x))^{-2}$. The first part of the following result is an analogue of \cite[Lemma 1.26]{BJK}.
\begin{lemma}\label{l:BJK-126}
Let  $f:(0,\infty)\to [0, \infty)$  be bounded on $(0,M]$ for every $M>0$, and $\lim_{t\to \infty}f(t)/t=0$. 

\noindent
(a) Let $D$ be a bounded open set,  $h>0$  a locally bounded function on $D$ such that $h\to \infty$ at $\partial D$ and
\begin{equation}\label{e:BJK-126}
\sup_{x\in D}\int_D G_D(x,y)h(y)\, dy <\infty.
\end{equation}
Then $f\circ h\in \JJ$. 

\noindent (b) Let $D$ be a bounded $C^{1,1}$ open set. Then $x\mapsto f(V(\d(x))^{-2})$ is in the Kato class $\JJ$.

\noindent (c)  Let $D$ be a bounded $C^{1,1}$ open set and let  $U:(0,\infty)\to [0,\infty)$ satisfy condition (U4). If
\begin{equation}\label{e:Kato-U}
\lim_{s\to 0} U(s) V(s)^2=0,
\end{equation}
then $x\mapsto U(\d(x))$ is in the Kato class $\JJ$.
\end{lemma}
\proof
(a) We will take advantage of the equivalence of (i) and (ii)  of \cite[Theorem 16.8]{Sch}. Denote $c:=\sup_{x\in D}\int_D G_D(x,y)h(y)\, dy$ and let $\eta>0$. There is $t_0>0$ such that $f(t)/t<\frac{\eta}{c}$ for every $t\ge t_0$. Also, since $h\to \infty$ at $\partial D$ there is $F\subset\subset D$ such that $h>t_0$ on $D\setminus F$ and since $h$ is locally bounded we have $M:=\sup_F h<\infty$. Hence
\begin{align*}
	\sup_{x\in D}\int_D G_D (x,y)f(h(y))dy&\le\sup_{x\in D}\int_F G_D (x,y)f(h(y))dy+\sup_{x\in D}\int_{D\setminus F} G_D (x,y)f(h(y))dy\\
	&\le ({\textstyle\sup_{(0,M]}}f)\, \sup_{x\in D}\mathbb{E}_x[\tau_D]+\eta<\infty,
\end{align*}
i.e. we have property $(a)$ of $(ii)$ in \cite[Theorem 16.8]{Sch}. Note that $\mathbf{1}\in \JJ$ since $D$ is bounded so there is $w_\eta\in L_+^1(D)$ and $\delta>0$ such that for all $B\subset D$ with $\int_Bw_\eta<\delta$ we have $\sup_{x\in D}\int_{B} G_D (x,y)dy<\frac{\eta}{\sup_{(0,M]}f}$. Hence, for all such $B$ it holds that
\begin{align*}
	\sup_{x\in D}\int\limits_B G_D (x,y)f(h(y))dy&\le \sup_{x\in D}\int\limits_{B\cap F} G_D (x,y)f(h(y))dy+\sup_{x\in D}\int\limits_{B\setminus F} G_D (x,y)f(h(y))dy	\\
	&\le ({\textstyle\sup_{(0,M]}}f)\left(\sup_{x\in D}\int\limits_{B} G_D (x,y)dy\right)+\eta\le 2\eta.
\end{align*}
Since $\eta$ was arbitrary we have  $(b)$ of $(ii)$ in \cite[Theorem 16.8.]{Sch}, i.e. $f\circ h\in\JJ$.

\noindent (b) This follows immediately from (a) by using \eqref{e:G-kappa-le-1} and \eqref{e:killing-fn-estimate}.

\noindent (c) Define $f(t):=U(V^{-1}(t^{-1/2}))$ so that $f(V(t)^{-2})=U(t)$. By the assumption on $U$, the function $f$ is locally bounded. Moreover, by using the substitution  $t=V(s)^{-2}$ and the assumption \eqref{e:Kato-U}, we get
$$
\lim_{t\to \infty}\frac{f(t)}{t}=\lim_{s\to 0} \frac{f( V(s)^{-2})}{V(s)^{-2}}=\lim_{s\to 0}U(s)V(s)^2=0.
$$
The claim now follows from (b).
\qed

\subsection{Generalized normal derivative, modified Martin kernel and equivalent formulation of the weak dual solution}\label{ss:gnd}
We now invoke the powerful recent result from \cite{KKLL} on boundary regularity of the solution of the equation
$$
\begin{array}{rcll}
  -Lu(x)&=&\psi(x)& \quad \text{in } D\\
   u&=&0 &\quad \text{in }D^c
\end{array}
$$
where $\psi$ is a bounded continuous function on $D$. It is proved in \cite[Theorem 1.2]{KKLL} (see also \cite[Theorem 3.10]{KKLL}), that $u=G_D \psi$ is the (viscosity) solution of the above equation, $u/V(\delta_D)\in C^{\gamma}(D)$, and
$$
\left \|\frac{u}{V(\delta_D)}\right\|_{C^{\gamma}(D)} \le C \|\psi\|_{\infty},
$$
for some constants $\gamma>0$ and $C>0$ depending only on $d$, $D$ and $\phi$. Here $C^{\gamma}(D)$ is the space of $\gamma$-H\"older continuous functions on $D$ with the corresponding H\"older norm.  It follows that $u/V(\delta_D)$ can be continuously extended to $\overline{D}$. In particular, for any bounded and continuous function $\psi:D\to \R$ and for every $z\in \partial D$, there exists a finite limit
\begin{equation}\label{e:normal-derivative}
\frac{d}{dV}(G_D \psi)(z):=\lim_{y\to z, y\in D}\frac{G_D \psi(y)}{V(\delta_D(y))}.
\end{equation}
We can think of $d(G_D \psi)/dV$ as the generalized normal derivative of the function $G_D \psi$ -- instead of the distance function $\delta_D$ we use $V(\delta_D)$.

If $\psi$ is nonnegative and has compact support, then $G_D \psi$ is regular harmonic in $D\setminus \mathrm{supp}(\psi)$. By \cite[Theorem 1.1]{KSV16}, for any $x\in D$, there exists a finite limit 
$$
\lim_{y\to z, y\in D} \frac{G_D \psi(y)}{G_D(x,y)}.
$$
Combining with \eqref{e:normal-derivative}, we see that for every $x\in D$ and every $z\in \partial D$, there exists
\begin{equation}\label{e:modified-martin}
K_D(x,z):=\lim_{y\to z, y\in D}\frac{G_D(x,y)}{V(\delta_D(y))}.
\end{equation}
We call $K_D(x,z)$ a modified Martin kernel, because given $x_0\in D$, we have that
\begin{equation}\label{e:K-and-M}
\frac{K_D(x,z)}{K_D(x_0,z)}=\lim_{y\to z, y\in D}\frac{\frac{G_D(x,y)}{V(\delta_D(y))}}{\frac{G_D(x_0,y)}{V(\delta_D(y))}}   =\lim_{y\to z}\frac{G_D(x,y)}{G_D(x_0,y)}=M_D(x,z).
\end{equation}

\begin{lemma}\label{l:derivative-alt}
 Let $D$ be a bounded open set and let $\psi:D\to \R$ be a bounded function with compact support and set $u=G_D \psi$. Then
$$
\frac{d}{dV}u(z)=\int_D K_D(y,z)\psi(y)\, dy. 
$$
\end{lemma}
\proof  Let $2\epsilon=\mathrm{dist}(\mathrm{supp}(\psi), \partial D)$,  $z\in \partial D$, and $x\in D$ such that $|x-z|<\epsilon$.
By using \eqref{e:G_D-estimate}, we get that for $y\in \mathrm{supp}(\psi)$, 
$$
\frac{G_D(x,y)}{V(\delta_D(x))}\le c \frac{ V(|x-y|)}{|x-y|^d}\le c\frac{V(\mathrm{diam}(D))}{\epsilon^d}.
$$
Thus we can use the bounded convergence theorem to conclude from \eqref{e:modified-martin}   that
$$
\frac{d}{dV}u(z)=\lim_{x\to z, x\in D}\frac{G_D \psi(x)}{V(\delta_D(x))}=\lim_{x\to z, x\in D}\int_D \frac{G_D(x,y)}{V(\delta_D( x ) )}\psi(y)\, dy = \int_D K_D(y,z)\psi(y)\, dy.
$$
\qed

Recall the weak dual formulation \eqref{e:weak-dual} of the semilinear problem \eqref{e:4}. We will now rewrite the last two integrals in \eqref{e:weak-dual}. Let $\psi\in C_c^{\infty}(D)$ and set $\varphi=G_D \psi$. First, by using \eqref{e:Lu} we see that
$$
\int_{ D^c}\int_{D}P_D(x,z)\psi(x)\, dx\, \lambda(dz)=-\int_{ D^c} (-L\varphi(z))\, \lambda(dz).
$$ 
Further, for $\mu \in \MM(\partial D)$, let $\wt{\mu}(dz):=K_D(x_0,z)\mu(dz)$. By Lemma \ref{l:derivative-alt} and \eqref{e:K-and-M}
\begin{align*}
\int_{\partial D}\int_{D}M_D(x,z)\psi(x)\, dx\,\wt{\mu}(dz)=\int_{\partial D}\int_{D}K_D(x,z)\psi(x)\, dx\, \mu(dz)=\int_{\partial D}\frac{d}{dV}\varphi(z)\mu(dz)\, .
\end{align*}
Since $\psi=-L\varphi$, we see that the function $u$ is a weak dual solution of the problem \eqref{e:4} if and only if
$$
\int_D u(x)(-L\varphi)(x)\, dx=\int_D f(x,u(x))\varphi(x)\, dx -\int_{\overline{D}^c} (-L\varphi(z))\, \lambda(dz) +\int_{\partial D}\frac{d}{dV}\varphi(z){\mu}(dz)\, .
$$
This formulation of a solution to the problem \eqref{e:4} in bounded $C^{1,1}$ open sets can be found in  \cite{Aba15a}  in the case of the fractional Laplacian.

\subsection{Another boundary operator}\label{ss:boundary-operator}
Following  \cite[Subsection 1.2]{Aba15a} (see also \cite[(2.2), Appendix B]{CGCV20})  we now introduce another boundary operator. For a measure $\mu \in \MM(\partial D)$ set $K_D\mu(x):=\int_{\partial D} K_D(x,z)\mu(dz)$, $x\in D$. 
Note that by Remark \ref{r:elementary-lemma}(i), $K_D(x_0, \cdot)$ is continuous on $\partial D$. In the context of the Proposition \ref{p:boundary-limit}, let $\mu(d\zeta):=f(\zeta)\sigma(d\zeta)$, $\wt{\mu}(d\zeta):=K_D(x_0,\zeta)\mu(d\zeta)$ and $\nu(d\zeta):=K_D(x_0,\zeta)\sigma(d\zeta)$. Then

\begin{equation}\label{e:boundary-limit}
\lim_{x\to z, x\in D}\frac{K_D \mu(x)}{K_D \sigma(x)}=\lim_{x\to z, x\in D}\frac{M_D \wt{\mu}(x)}{M_D \nu(x)}=\frac{K_D (x_0,z)f(z)}{K_D (x_0,z)}=f(z).
\end{equation}

For $u:D\to \R$ and $z\in \partial D$, let
$$
E_D u(z):=\lim_{x\to z, x\in D} \frac{u(x)}{K_D \sigma (x)},
$$
whenever the limit exists and is finite. 

\begin{remark}\label{r:elementary-lemma}
{\rm We will need the following elementary calculations several times below.
\begin{enumerate}[leftmargin=0cm,itemindent=.5cm,labelwidth=\itemindent,labelsep=0.03cm,align=left, label=(\roman*)]
\item Let $u:D\to \R$ be a function and assume that for every $z\in \partial D$ there exists a finite limit
\begin{equation}\label{e:elementary-lemma}
\wt{u}(z):=\lim_{x\to z, x\in D} u(x).
\end{equation}
Then, by applying the usual $2\varepsilon$-argument, it follows that $\wt{u}:\partial D\to \R$ is continuous.
\item Assume further that $D$ is bounded and $\wt{u}(z)=0$ for all $z\in \partial D$. Then convergence in \eqref{e:elementary-lemma} is uniform in the sense that for every $\epsilon >0$ there exists a compact set $F\subset D$ such that $|u(x)|<\epsilon $ for all $x\in D\setminus F$. Indeed, due to compactness of $\partial D$ we easily find a finite cover $V:=\cup_{i=1}^n B(z_i,r_i)$, $z_i\in\partial D$, of $\partial D$ such that $|u|\le \varepsilon$ on $D\cap V$.
\end{enumerate}
}
\end{remark}

\begin{prop}\label{p:boundary-operators}
Let $u:D\to \R$. If $E_D u(z)$ exists for every $z\in\partial D$, then $W_D u$ exists and 
\[
W_D u(dz)=E_D u(z) K_D(x_0, z)\sigma(dz).
\]
\end{prop}
\proof Assume that $E_D u(z)$ exists for every $z\in \partial D$. By Remark \ref{r:elementary-lemma}(i), $E_D u$ is continuous on $\partial D$.
Let $\nu(dz)=K_D(x_0,z)\sigma(dz)$, $\mu(dz)=E_Du (z)\nu(dz)$ and
$$
v(x):=M_D \mu(x)=\int_{\partial D}M_D(x,z)E_D u(z)\nu(dz)=\int_{\partial D} K_D(x,z) E_D u(z) \sigma(dz).
$$
By \eqref{e:boundary-limit}, for every $z\in \partial D$,
$$
\lim_{x\to z, x\in D}\frac{v(x)}{K_D \sigma(x)}=E_D u(z),
$$
hence $E_D v= E_D u$, so that $\lim_{x\to z, x\in D}(u(x)-v(x))/K_D\sigma(x)=0$ for every $z\in \partial D$. By Remark \ref{r:elementary-lemma}(ii), this implies that for every $\epsilon>0$ there exists a compact set $F\subset D$, such that
$$
\frac{|u(x)-v(x)|}{K_D\sigma (x)}<\epsilon, \quad \textrm{for all }x\in D\setminus F.
$$
Since $K_D\sigma$ is a nonnegative harmonic function, the same proof as \cite[Lemma 1.16]{BJK} gives that $W_D(u-v)=0$. Notice that the set of functions on $D$ for which $W_D $ is defined is a vector space and  $W_D $ is linear on that space.  We conclude that $W_D u$ exists and $W_D u =W_D v+W_D(u-v)=W_D v = W_D(M_D\mu)=\mu$  by \eqref{e:W_D(G_D)}.
\qed

\section{The semilinear problem in bounded $C^{1,1}$ open set}\label{s:semilinear-C11}
\subsection{Corollary \ref{c:bogdan} revisited}

Recall that in Corollary \ref{c:bogdan} we assumed that the function  $f:D\times \R\to \R$ satisfies \textbf{(F)} with $\Lambda$ nondecreasing and that $G_D\rho\in C_0(D)$ and $G_D(\rho \Lambda(2\overline{g}))\in C_0(D)$, where $\overline{g}=P_D|\lambda|+M_D|\mu|$. We give sufficient conditions for these assumptions in case of a bounded $C^{1,1}$ open set. We will additionally assume that $\rho(x)=W(\d(x))$ for a function $W:(0,\infty)\to [0,\infty)$ and that $\Lambda$ satisfies the following doubling condition: There exists $C\ge 1$ such that 
\begin{equation}\label{e:doubling-Lambda}
\Lambda (2t)\le C \Lambda(t), \quad t>0.
\end{equation}
This implies that for all $c_1>1$ there exists $c_2=c_2(C, c_1)$ such that
$$
\Lambda (c_1 t)\le c_2\Lambda (t), \quad \quad t>0,
$$
which can be rewritten as follows: For every $\wt{c}_1\in (0,1)$, there exists $\wt{c}_2>0$ such that
\begin{equation}\label{e:doubling-Lambda-reverse}
\Lambda (\wt{c}_1t)\ge \wt{c}_2\Lambda (t), \quad t>0.
\end{equation}
Secondly, assume that
$$
\overline{g}(x)\preceq \frac{V(\d(x))}{\d(x)}, \quad x\in D.
$$
 By \eqref{e:P_Dg2} and \eqref{e:upper-bound-M}, this will be the case provided $\mu(dz)=h(z)\sigma(dz)$ for a continuous function $h:\partial D\to \R$, and $\lambda(dy)=g(y)dy$ with $|g(y)|\preceq \wt{U}(\dc(y))$ where  $\wt{U}$ is nonnegative, bounded on compact subsets of $(0,\infty)$ and satisfies  \eqref{e:poisson-finite}.  Then we have 
$$
\rho(x)\Lambda(2\overline{g})(x)\le c W(\d(x))\Lambda\left(\frac{V(\delta_D(x))}{\delta_D(x)}\right), \quad x\in D,
$$
for some $c>0$. By using Lemma \ref{l:BJK-126}(c), we see that $G_D(\rho \Lambda(2\overline{g}))\in C_0(D)$ if
$$
\lim_{t\to 0} W(t)\Lambda\left(\frac{V(t)}{t}\right)V(t)^2=0\, ,
$$
while $G_D\rho \in C_0(D)$ if $\lim_{t\to 0}W(t)V(t)^2=0$. 

In the case of the fractional Laplacian, $W(t)=t^{-\beta}$ and $\Lambda(t)=t^p$, these two conditions become $\beta+p(1-\alpha/2)<\alpha$. 

\subsection{Theorem \ref{t:nonpositive-f} in bounded $C^{1,1}$ open set}
In this subsection we revisit Theorem \ref{t:nonpositive-f}(ii)  in case of a bounded $C^{1,1}$ open set $D$. Recall that the assumptions of that theorem were that $f:D\times \R\to (-\infty, 0]$ satisfies \textbf{(F)} with $G_D\rho \in C_0(D)$, $f(x,0)=0$ and the function $\Lambda $ is nondecreasing.
As in the previous subsection, we will additionally assume that $\rho(x)=W(\d(x))$ for a function $W:(0,\infty)\to [0,\infty)$ and that $\Lambda$ satisfies  the doubling condition \eqref{e:doubling-Lambda}.

\begin{prop}\label{p:nonposiitve-f-revisited}
Let $D\subset \R^d$ be a bounded $C^{1,1}$ open set. Let $f:D\times \R\to (-\infty, 0]$ be a function that satisfies \textbf{(F)} with $\rho(x)=W(\d(x))$, where $W:(0,\infty)\to [0, \infty)$ is bounded away from zero, and such that $\Lambda$ is a nondecreasing function  satisfying  the doubling condition \eqref{e:doubling-Lambda}. Assume that
\begin{equation}\label{e:W-Kato}
\lim_{t\to 0} W(t)V(t)^2=0\, .
\end{equation}
Let $\lambda (dy)=\wt{U}(\dc(y))dy$ where $\wt{U}:(0,\infty)\to [0,\infty)$ is bounded  on every compact subset of $(0,\infty)$  and satisfies \eqref{e:poisson-finite}, and let $\mu(dz)=h(z)\sigma(dz)$ where $h:\partial D\to[0,\infty)$ is continuous and not identically equal to zero. If for some $\eta>0$
\begin{equation}\label{e:int-criterion}
	\int_{0}^\eta W(t)V(t)\Lambda\left(\frac{V(t)}{t}\right)dt<\infty,
\end{equation}
then the semilinear problem \eqref{e:4} has a nonnegative weak dual solution $u\in L^1(D)\cap C(D)$.
\end{prop}
\proof We first note that the assumption  \eqref{e:W-Kato} implies by Lemma \ref{l:BJK-126}(c)  that $\rho=W(\d)\in \JJ$, and thus by Subsection \ref{ss:aux}, $\rho\in C_0(D)$. Hence, in order to see that the semilinear problem \eqref{e:4} has a nonnegative solution it suffices to check that $\rho \Lambda(g)G_D \1\in L^1(D)$ where $g=P_D\lambda +M_D\mu$. By \eqref{e:P_Dg2} and \eqref{e:upper-bound-M} there exists a constant $c_1>0$ such that 
$$
g(x)\le  c_1 \frac{V(\d(x))}{\d(x)}, \quad x\in D.
$$
Together with \eqref{e:doubling-Lambda} this implies that  
$$
\Lambda(g(x))\le \Lambda\left( c_1  \frac{V(\d(x))}{\d(x)}\right) \le  c_2  \Lambda\left(\frac{V(\d(x))}{\d(x)}\right), \quad x\in D,
$$
for some $ c_2 >0$.  Therefore, there is $ c_3>0$ such that
\begin{equation}\label{e:nonposiitve-f-revisited1}
\rho(x)\Lambda(g(x))G_D \1(x)\le  c_3  W(\delta_D(x)) \Lambda\left(\frac{V(\delta_D(x))}{\delta_D(x)}\right) V(\delta_D(x)), \quad x\in D.
\end{equation}
By using  boundedness of $W(\delta_D) \Lambda\left(\frac{V(\delta_D)}{\delta_D}\right) V(\delta_D)$ inside $D$ and the co-area formula near the boundary of $D$ with the assumption \eqref{e:int-criterion}  we see that
$$
\int_D W(\delta_D(x)) \Lambda\left(\frac{V(\delta_D(x))}{\delta_D(x)}\right) V(\delta_D(x))dx <\infty.
$$
Now it follows from \eqref{e:nonposiitve-f-revisited1} that $\rho \Lambda(g)G_D \1\in L^1(D)$.
\qed

\begin{remark}\label{r:nonposiitve-f-revisited}
{\rm
(a) Proposition \ref{p:nonposiitve-f-revisited} allows a partial converse. Assume that  $f(x,t)=  - W(\delta_D(x))\Lambda( |t|)$ where $\Lambda:(0, \infty)\to (0, \infty)$ is a nondecreasing and unbounded function satisfying \eqref{e:doubling-Lambda}
Assume further that there exists $\eta_0>0$ such that for all $\eta\in(0,\eta_0]$
	\begin{equation}\label{e:int-criterion-2}
	\int_{0}^\eta W(t)V(t)\Lambda\left(\frac{V(t)}{t}\right)dt=+\infty.
	\end{equation}
Let  $\mu(d\zeta)=h(\zeta) \sigma(d\zeta)$ with nonnegative continuous $h$, $h\neq 0$. Then the semilinear problem \eqref{e:4}
does not have a nonnegative weak dual solution $u\in L^1(D)$ such that $E_D u$ is well defined. To show this, suppose that there exists  a nonnegative $u$ that solves  \eqref{e:4}. Then $u(x)=G_Df_u(x)+P_D\lambda (x)+M_D\mu(x)$ a.e. 
Since by assumption, $E_D u$ exists, by Proposition 4.8 $W_D u$ also exists and $W_D u(d\zeta)=E_D u(\zeta)K_D(x_0,\zeta)\sigma(d\zeta)$. On the other hand, since $u=G_D f_u+P_D\lambda +M_D\mu$, we have by (2.15) that $W_D u=W_D(M_D\mu)=\mu$. Since $\mu(d\zeta)=h(\zeta)\sigma(d\zeta)$, we get 
$$
E_D u(\zeta)=\frac{h(\zeta)}{K_D(x_0, \zeta)} \quad \sigma(d\zeta)- \textrm{a.e.}
$$
Choose $z\in \partial D$ such that $E_D u(z)=h(z)/K_D(x_0,z)>0$. Since
$$
E_D u(z)=\lim_{x\to z, x\in D}\frac{u(x)}{K_D \sigma(x)},
$$
there exists $\epsilon>0$ such that
$$
u(x)\ge \frac{1}{2}E_D u( z ) K_D\sigma(x)=\frac{1}{2}\frac{h(z)}{K_D(x_0,z)}K_D\sigma(x)=  c_1 K_D\sigma(x), \qquad \textrm{ for all }x\in D\cap B(z, \epsilon),
$$
where $c_1=c_1(z,h)>0$. By using \eqref{e:G_D-estimate} and \eqref{e:modified-martin} to get the same estimate of $K_D(x,z)$ as the one of $M_D(x,z)$ in \eqref{e:M_D-estimate}, we see in the same way as for \eqref{e:M_D-sigma} that
$$
K_D\sigma(x)\asymp \frac{V(\delta_D(x))}{\delta_D(x)}, \quad x\in D.
$$
This implies that there exists $c_2=c_2(z,h)>0$ such that
$$
u(x)\ge c_2\frac{V(\delta_D(x))}{\delta_D(x)},  \qquad \textrm{ for all }x\in D\cap B(z, \epsilon).
$$
Therefore, by using \eqref{e:doubling-Lambda-reverse}  this implies that for some $c_3>0$
$$
\Lambda(u(y))\ge \Lambda\left( c_2\frac{V(\delta_D(y))}{\delta_D(y)}\right)   c_3\Lambda\left(\frac{V(\delta_D(y))}{\delta_D(y)}\right), \quad \textrm{for all } y\in D\cap B(z, \epsilon).
$$
Choose $x\in D$ so that $\delta_D(x)\asymp |x-y|\asymp 1$ whenever $y\in D\cap B(z,\epsilon)$. By \eqref{e:G_D-estimate}, there exists $c_4>0$ such that $G_D(x,y)\ge c_4 V(\delta_D(y))$. Hence,
$$
G_D f_u(x)=\int_D G_D(x,y) f(y, u(y))dy  \le -    c_3 c_4\int_{D\cap B(z, \epsilon)} V(\delta_D(y)) W(\delta_D(y))\Lambda\left(\frac{V(\delta_D(y))}{\delta_D(y)}\right)dy.
$$
By use of the co-area formula it follows that the last integral is equal to some constant multiplied by
$$
\int_0^{\epsilon}V(t)W(t)\Lambda\left(\frac{V(t)}{t}\right)dt.
$$
By \eqref{e:int-criterion-2}  it follows that $G_D f_u(x)= -\infty $ for points $x$ in some open subset of $D$. This is a contradiction with $G_D f_u   > -\infty $ a.e. which follows from  $u\ge 0$, $P_D \lambda<\infty$ and $M_D \mu<\infty$.

(b) Note that the power function $\Lambda(t)=t^p$ is increasing and satisfies the doubling condition \eqref{e:doubling-Lambda}. Assume that $W(t)=t^{-\beta}$ and the underlying process is an isotropic $\alpha$-stable process (so that $V(t)= t^{\alpha/2} $). Then \eqref{e:W-Kato} reads  $\beta < \alpha$, while the integral criterion \eqref{e:int-criterion} is equivalent to $\beta+p(1-\alpha/2)<1+\alpha/2$. 
In case $f(x,t)= -  t^p $, we see that the problem \eqref{e:4} has a nonnegative solution $u$ if $p<(2+\alpha)/(2-\alpha)$, while in case $p\ge (2+\alpha)/(2-\alpha)$  a nonnegative solution $u$ such that $E_D u$ is well defined does not exist.
}
\end{remark}

\subsection{Extending Corollary \ref{c:bogdan} to a wider class of nonnegative nonlinearities}

Our next goal is to extend the results of Corollary \ref{c:bogdan} for nonnegative nonlinearities $f$. Unlike Theorem \ref{t:nonpositive-f}, this approach relies heavily on the estimates of Green and Poisson potentials in bounded $C^{1,1}$ domains. 

\begin{thm}\label{thm:superlinear}
Let $f:D\times\mathbb R\to\mathbb [0,\infty)$ be a function, nondecreasing in the second variable, satisfying \textbf{(F)}, with $\rho=W(\d)$ for some function $W:(0,\infty)\to[0,\infty)$, $\Lambda$ nondecreasing and satisfying  the doubling condition \eqref{e:doubling-Lambda}.  Let $\lambda\in \mathcal M(D^c)$ be a nonnegative measure which is absolutely continuous with respect to the Lebesgue measure with density $\widetilde U(\delta_{D^c})$, where $\widetilde U:(0,\infty)\to[0,\infty)$ is a  function bounded on compact subsets of $(0,\infty)$ satisfying \eqref{e:poisson-finite}. Let $h:\partial D\to[0,\infty)$ be a continuous function and let $\mu(d\zeta)=h(\zeta) \sigma(d\zeta)$ be a measure on $\partial D$. Suppose that one of the following conditions  hold: 
\begin{enumerate}[label=(\roman*)]
\item the function $t\mapsto W(t)\Lambda\left(\frac{V(t)}{t}\right)$, $t>0$, satisfies the conditions \textbf{(U)};
\item $h\equiv 0$ and the function $W\Lambda(\widetilde U)$ satisfies the conditions \textbf{(U)}. Morover assume that
\begin{align}
&\int_0^{\mathrm{diam}(D)}\frac{\widetilde U(t)}{V(t)(s+t)}\,dt\preceq \frac{\widetilde U(s)}{V(s)},\label{e:poisson-potential-dominated}\\
\begin{split}\label{e:green-potential-dominated}
&\int_{0}^s W(t)V(t)\Lambda(\widetilde U(t))dt\preceq \frac{s\widetilde U(s)}{V(s)},\\
&\int_{s}^{\mathrm{diam}(D)} \frac{W(t)V(t)\Lambda(\widetilde U(t))}t\, dt\preceq \frac{\widetilde U(s)}{V(s)},
\end{split}
\end{align}
where the constants do not depend on $0<s\leq \tfrac{\mathrm{diam}(D)}2$.
\end{enumerate}
Then there exists a constant $m_1>0$ such that for every $m\in[0,m_1]$ the semilinear problem
\begin{equation}\label{e:superlinear}
\begin{array}{rcll}
		- Lu(x)&=& m f(x,u(x))& \quad \text{in } D\\
		u&=&\lambda& \quad \text{in } \overline{D}^c\\
		W_Du&=&\mu&\quad \text{on }\partial D
\end{array}
\end{equation}
has a nonnegative weak dual solution $u\in L^1(D)$.
\end{thm}
\proof
First we prove the theorem under assumption $(i)$. Since $f$ is nonnegative, the function $u_0=P_D\lambda+M_D\mu$ is a subsolution to \eqref{e:superlinear}. Recall from \eqref{e:P_Dg2} and \eqref{e:upper-bound-M} that there exists a constant $ c_1>0$ such that
\[
u_0(x)\le  c_1\frac{V(\d(x))}{\d(x)}, \ x\in D.
\]
Next we construct a supersolution $\overline{u}$ for \eqref{e:superlinear} of the form 
\[
\overline{u}(x)= c_2\frac{V(\d(x))}{\d(x)},
\] 
i.e. find a constant $ c_2>c_1$ such that 
\begin{equation}\label{e:superlinear-supersolution}
\overline{u}(x)\ge mG_Df_{\overline{u}}(x) + u_0(x),\ x\in D,
\end{equation}
for $m$ small enough. To be exact, we show that for every $ c_2> c_1$ there exists $m_1>0$ such that \eqref{e:superlinear-supersolution} holds for every $m\in [0, m_1 ]$. Fix $ c_2> c_1$. First note that by {\bf (F)} and the doubling property \eqref{e:doubling-Lambda} for $\Lambda$ we have 
\[
f\left(x, c_2\frac{V(\d(x))}{\d(x)}\right)\le W(\d(x))\Lambda\left( c_2\frac{V(\d(x))}{\d(x)}\right)\leq c_3 W(\d(x))\Lambda\left(\frac{V(\d(x))}{\d(x)}\right)
\]
for some constant $c_3>0$. Now by Proposition \ref{p:green-potential-estimate} there exists $c_4>0$ such that
\[
G_Df_{\overline{u}}(x)\leq  c_3 G_D\left[W(\d)\Lambda\left(\frac{V(\d)}{\d}\right)\right](x)\leq  c_4 \frac{V(\d(x))}{\d(x)}. 
\]
By choosing $m_1= \frac{c_2-c_1}{c_4}$ we get that for every $m \le  m_1$
\[
mG_Df_{\overline{u}}(x)+u_0(x)\leq ( m c_4+c_1) \frac{V(\d(x))}{\d(x)}=\frac{ m c_4+c_1}{ c_2}\overline{u}(x)\leq \overline{u}(x). 
\]
Now we can apply the classical iteration scheme in the following way: For $k\in\mathbb N$ let $u_k$ be the weak $L^1$ solution to the linear problem 
\begin{equation*}
\begin{array}{rcll}
		- Lu_k(x)&=& m f(x,u_{k-1}(x))& \quad \text{in } D\\
		u_k&=&\lambda& \quad \text{in } \overline{D}^c\\
		W_Du_k&=&\mu&\quad \text{on }\partial D.
\end{array}
\end{equation*}

The constructed sequence $(u_k)_k$ is nondecreasing and dominated by $\overline{u}$. To see this, take $x\in D$. Since $f$ is nonnegative, we have that 
\[
u_1(x)-u_0(x)=mG_Df_{u_0}(x)\geq 0.
\]
Furthermore, since $f$ is nondecreasing in the second variable and $u_0\le \overline{u}$, we have that
\[
u_1(x)=mG_Df_{u_0}(x)+u_0(x)\leq mG_Df_{\overline{u}}(x)+u_0(x)\leq \overline{u}(x). 
\]
Assume now that $u_{k-1}(x)\leq u_k(x)\leq \overline{u}(x)$ for some $k\in\mathbb N$. This implies that $f_{u_{k-1}}(x)\leq f_{u_{k}}(x)\leq f_{\overline{u}}(x)$, so
\[
u_{k+1}(x)-u_{k}(x)=mG_Df_{u_k}(x)-mG_Df_{u_{k-1}}(x)\geq 0
\]
and
\[
u_{k+1}(x)= mG_Df_{u_k}(x)+u_{0}(x)\le mG_Df_{\overline u}(x)+u_{0}(x)\le \overline{u}(x). 
\]
The claim now follows by induction.

Therefore, we can define a pointwise limit $u:=\lim_{k\to\infty} u_k$ which, by the monotone convergence theorem and the continuity of $f$ in the second variable, satisfies
\begin{align*}
u(x)&=\lim_{k\to\infty} \int_D f(y,u_{k}(y))G_D(x,y) dy + u_0(x)\\
&= \int_D \lim_{k\to\infty}f(y,u_{k}(y))G_D(x,y) dy + u_0(x)\\
&= \int_D f(y,u(y))G_D(x,y) dy + u_0(x),
\end{align*}
i.e. $u$ is a weak $L^1$ solution to \eqref{e:superlinear}.

Next, consider the proof of the theorem under the assumptions $(ii)$. Note that we only need to find a supersolution $\overline{u}\geq u_0=P_D\lambda$ satisfying \eqref{e:superlinear-supersolution}. The rest of the proof then follows from the proof of $(i)$. Note first that \eqref{e:P_Dg} and \eqref{e:poisson-potential-dominated} imply that there exists a constant $ c_5  >0$ such that 
\[
u_0(x)\leq  c_5  \widetilde{U}(\delta_D(x)),\ x\in D.
\]
Therefore, in this case we fix a constant $c_6>c_5$ and show that the function $\overline{u}$ of the form
\[
\overline{u}(x)=c_6 \widetilde U(\delta_{D}(x)),\ 
\]
is indeed a supersolution  to \eqref{e:superlinear} for $m$ small enough. As in the previous case, by {\bf (F)} and the doubling property for $\Lambda$ 
\[
f\left(x, c_6 \widetilde U(\delta_{D}(x))\right)\leq W(\d(x))\Lambda\left(c_6 \widetilde U(\delta_{D}(x))\right)\leq c_7  W(\d(x))\Lambda\left(\widetilde U(\delta_{D}(x))\right)
\]
for some constant $c_7>0$. Now by Proposition \ref{p:green-potential-estimate} and \eqref{e:green-potential-dominated} it follows that 
\[
G_Df_{\overline{u}}(x)\leq  c_7  G_D\left[W(\d)\Lambda\left(\widetilde U(\delta_{D}(x))\right)\right](x)\leq  c_8  \widetilde U(\delta_{D}(x)). 
\]
By choosing $m_1=\frac{c_6-c_5}{c_8}$ we get that for every $ m\le m_1$
\[
mG_Df_{\overline{u}}(x)+u_0(x)\leq ( mc_8+c_5) \widetilde U(\delta_{D}(x))=\frac{mc_8+c_5}{c_6}\overline{u}(x)\leq \overline{u}(x). 
\]
\qed

Assume that functions $W$ and $\Lambda$ satisfy \eqref{e:int-criterion}, $W$ satisfies conditions (U2)-(U4), and $\Lambda$  is nondecreasing  and  satisfies the doubling condition \eqref{e:doubling-Lambda}.
Then the function $U(t)=W(t)\Lambda\left(\frac{V(t)}{t}\right)$ satisfies conditions \textbf{(U)}. Indeed, since $W$ is almost nonincreasing and $\Lambda$ is nondecreasing it follows that
\begin{align*}
\frac{W(t)\Lambda\left( \frac{V(t)}{t}\right)}{W(s)\Lambda\left( \frac{V(s)}{s}\right)}\preceq \frac{\Lambda\left( \frac{V(t)}{t}\right)}{\Lambda\left( \frac{V(s)}{s}\right)}\overset{\eqref{e:wsc-V}}\preceq \frac{\Lambda\left(  \wt a_2\frac{V(s)}{s}\right)}{\Lambda\left( \frac{V(s)}{s}\right)}\overset{\eqref{e:doubling-Lambda}}{\preceq} 1,\ \ s< t\le 1.  
\end{align*}
Furthermore, since $W$ satisfies the reverse doubling condition \eqref{e:doubling-U} and $\Lambda$ is nondecreasing, we have that 
\begin{align*}
\frac{W(t)\Lambda\left( \frac{V(t)}{t}\right)}{W(2t)\Lambda\left( \frac{V(2t)}{2t}\right)}\preceq \frac{\Lambda\left( \frac{V(t)}{t}\right)}{\Lambda\left( \frac{V(2t)}{2t}\right)}\overset{\eqref{e:wsc-V}}\preceq \frac{\Lambda\left(\frac{V(t)}{t}\right)}{\Lambda\left(  \wt a_1 2^{\delta_1-1}\frac{V(t)}{t}\right)}\overset{\eqref{e:doubling-Lambda}}{\preceq} 1,\ \ t\in(0,1).  
\end{align*}
Finally, note that $U$ is bounded away from zero, since both $W$ and $t\mapsto \frac{V(t)}{t}$ satisfy (U4) and $\Lambda$ is nondecreasing. 
Similarly, note that the function $U=W\Lambda(\widetilde U)$ satisfies conditions (U2)-(U4) if we additionally assume that $\widetilde U$ satisfies (U2)-(U4).

\begin{remark}
{\rm
\begin{enumerate}[leftmargin=0cm,itemindent=.5cm,labelwidth=\itemindent,labelsep=0.03cm,align=left, label=(\roman*)]
\item Consider the isotropic $\alpha$-stable case and take $\Lambda(t)=t^p$ and $W(t)=t^{-\beta_1}$ for some $p>0$ and $\beta_1\ge0$, as in Remark \ref{r:nonposiitve-f-revisited}. The function $U(t)=W(t)\Lambda\left(\frac{V(t)}{t}\right)$ satisfies conditions \textbf{(U)} if and only if $\beta_1+p(1-\alpha/2)<1+\alpha/2$. Hence, if $f(x,t)=t^p$, then Theorem \ref{thm:superlinear} holds for $p<\frac{2+\alpha}{2-\alpha}$. 
\item When $\widetilde U(t)=t^{-\beta_2}$, the function $W\Lambda(\widetilde U)$ satisfies conditions \textbf{(U)} if and only if $\beta_1+p\beta_2<1+\alpha/2$. The condition \eqref{e:poisson-potential-dominated} is satisfied for $\beta_2<1-\alpha/2$. When $\beta_1=0$ the conditions in \eqref{e:green-potential-dominated} are satisfied when $\beta_2(p-1)\leq \alpha$. Since $\beta_2<1-\alpha/2$ we have that $\frac{\alpha}{\beta_2}+1<\frac{1+\alpha/2}{\beta_2}$, so Theorem \ref{thm:superlinear} states that the solution exists for  $p<\frac{\alpha}{\beta_2}+1$.
\end{enumerate}
}
\end{remark}

\section{Appendix}

\subsection{Approximation of excessive functions}
Let $(X_t, \P_x)$ be a Hunt process on a locally compact space $D$ and let $(P_t)_{t\ge 0}$ denote its semigroup. Let $U$ be the potential operator of $X$, that is
$$
Uf(x)=\E_x \int_0^{\infty} f(X_t)\, dt=\E_x\int_0^{\zeta}f(X_t)dt=\int_0^{\infty}P_t f(x)dt.
$$
Here $\zeta$ denotes the lifetime of the process. We assume that $X$ is transient in the sense that there exists a nonnegative measurable function $h$ such that $0<Uh<\infty$, see \cite[p.86]{CW}, and also that $(P_t)$ is strongly Feller. What follows essentially comes from \cite[Section 3.2]{CW}. Recall that a measurable function $f:E\to [0,\infty]$ is said to be excessive relative to $(P_t)_{t\ge 0}$ if $f\ge P_t f$ for all $t\ge 0$ and $f=\lim_{t\to 0}P_t f$ (see for example \cite[Section 2.1]{CW}). 
\begin{lemma}\label{l:1}  Suppose that $f$ is excessive, $P_t f <\infty$ for all $t\ge 0$ and $\lim_{t\to \infty}P_t f=0$. Then there exists a sequence $(g_n)_{n\ge 1}$ of nonnegative measurable functions such that $f=\uparrow \lim Ug_n$. Moreover, if $f$ is continuous and bounded, then one can choose $g_n$ to be continuous.
\end{lemma}
\proof This is proved as \cite[Theorem 6, p.82]{CW}. The function $g_n$ is given by
$$
g_n=n(f-P_{1/n}f).
$$
If $f$ is bounded, then $P_{1/n}f$ is continuous (by  the  strong Feller  property). If $f$ is also continuous, then $f-P_{1/n}f$ is continuous.\qed

\begin{remark}
{\rm 
Transience is not needed in this result. The assumption $P_t f <\infty$ is satisfied if $f<\infty$ since $P_t f\le f$. The assumption $\lim_{t\to \infty} P_t f=0$ is \emph{not} satisfied for harmonic functions (since they are invariant).
}
\end{remark}

\begin{prop}\label{p:approximation}
Let $f$ be excessive. If $(P_t)$ is transient, there exists a sequence $(g_n)_{n\ge 1}$ of bounded measurable functions such that $f=\uparrow \lim_{n\to \infty} U g_n$. Moreover, assume that there exists $h>0$ such that $0<Uh<\infty$ and $Uh$  is  continuous. If $f$ is continuous and $(P_t)$ is strongly Feller, then one can choose $g_n$ to be continuous.
\end{prop}
\proof Let $h_n=n h$ with $0<Uh<\infty$. and put 
$$
f_n=f\wedge Uh_n \wedge n.
$$
By \cite[Theorem 8, p.104]{CW}, $f_n$ is excessive (minimum of excessive function is excessive). Note that under additional assumptions, $f_n$ is continuous (and clearly bounded). By Lemma \ref{l:1}, there exists a sequence $(g_{nk})_{k\ge 1}$ such that $f_n=\uparrow \lim_{k\to \infty} U g_{nk}$. In fact,
$$
g_{nk}=k(f_n-P_{1/k}f_n) \le kn\, .
$$
Under additional assumptions, $g_{nk}$ are continuous. From the proof of Lemma \ref{l:1}, cf.~\cite[Theorem 6, p.82]{CW},
$$
U g_{nk}=k\int_0^{1/k}P_s f_n\, ds \le n\, .
$$
For each $n$, $Ug_{nk}$ increases with $k$ (this is part of Lemma \ref{l:1}); for each $k$, $U g_{nk}$ increases with $n$ (this follows from $f_n\le f_{n+1}$). Now, by \cite[Lemma 1, p.80]{CW},
$$
\uparrow \lim_{n\to \infty} f_n=\uparrow \lim_{n\to \infty} \uparrow \lim_{k\to \infty} U g_{nk}=\uparrow \lim_{n\to \infty} U g_{nn}.
$$
On the other hand, by the same \cite[Lemma 1, p.80]{CW} and monotone convergence
$$
\uparrow \lim_{n\to \infty} f_n=\uparrow \lim_{n\to \infty} \uparrow \lim_{t\downarrow 0} P_t f_n = \uparrow \lim_{t\downarrow 0} \uparrow \lim_{n\to \infty} P_t f_n \uparrow \lim_{t\downarrow 0}  P_t f=f.
$$ 
Therefore, by setting $g_n=g_{nn}$,
$$
f=\uparrow \lim_{n\to \infty} U g_n.
$$
\qed 

\subsection{Proofs of Propositions \ref{p:green-potential-estimate} and \ref{p:poisson-potential-estimate}}

Let $\epsilon >0$ be such that the map $\Phi:\partial D\times (-\epsilon, \epsilon)\to \R^d$ defined by $\Phi(y,\delta)=y+\delta\mathbf{n}(y)$ defines a diffeomorphism to its image, cf.~\cite[Remark 3.1]{AGCV19}. Here $\mathbf{n}$ denotes the unit interior normal.  Without loss of generality assume that $\epsilon<\mathrm{diam}(D)/20$.

\begin{lemma}\label{l:close-to-bdry}
Let $\eta<\epsilon$ and assume that  conditions \textbf{(U)} hold true. Then for any $x\in D$ such that $\delta_D(x)<\eta/2$, 
\begin{eqnarray}\label{e:close-to-bdry}
G_D\big(U(\delta_D)\1_{(\delta_D<\eta)}\big)(x)&\asymp &\frac{V(\d(x))}{\d(x)}\int_0^{\d(x)}U(t)V(t)\, dt + V(\d(x))\int_0^{\eta}U(t)V(t)\, dt \nonumber \\
& &+\ V(\d(x))\int_{3\d(x)/2}^{\eta} \frac{U(t)V(t)}{t}\, dt\, .
\end{eqnarray}
Further, $G_D\big(U(\delta_D)\1_{(\delta_D<\eta)}\big)(x)<\infty$ if and only if  the integrability condition \eqref{e:int-cond} holds true. 
\end{lemma}
\proof Let $r_0:=\mathrm{diam}(D)/10$. Fix $x\in D$ as in the statement and define
\begin{eqnarray*}
D_1&=& B(x, \delta_D(x)/2)\\
D_2&=&\{y: \delta_D(y)<\eta\}\setminus B(x, r_0)\\
D_3&=&\{y: \delta_D(y)<\delta_D(x)/2\}\cap B(x, r_0)\\
D_4&=& \{y:3\delta_D(x)/2 < \delta_D(y)<\eta\}\cap B(x,r_0)\\
D_5&=&\{y: \delta_D(x)/2<\delta_D(y)<3\delta_D(x)/2\}\cap (B(x, r_0)\setminus B(x, \delta_D(x)/2)).
\end{eqnarray*}
Thus we have that
$$
G_D\big(U(\delta_D)\1_{(\delta_D<\eta)}\big)(x)=\sum_{j=1}^5 \int_{D_j}G_D(x,y)U(\delta_D(y))\, dy =:\sum_{j=1}^5 I_j.
$$

\noindent
{\bf Estimate of $I_1$:} Under the almost nonincreasing condition \eqref{e:wus-U} and the doubling condition \eqref{e:doubling-U} it holds that
\begin{equation}\label{e:estimate-I1}
I_1 \preceq U(\delta_D(x))V(\delta_D(x))^2 \preceq \frac{V(\d(x))}{\d(x)}\int_0^{\d(x)}U(t)V(t)\, dt\, .
\end{equation}
Indeed, let $y\in D_1$. Then $\d(y)>\d(x)/2>|y-x|$ implying that 
$$
\left( 1\wedge \frac{V(\delta_D(x))}{V(|x-y|)}\right) \left( 1\wedge \frac{V(\delta_D(y))}{V(|x-y|)}\right) \asymp 1.
$$
Further, by using first \eqref{e:wus-U} and then \eqref{e:doubling-U} we have that
\begin{equation}\label{e:estimate-I1b}
U(\d(y)) \le  c_1   U(\d(x)/2)\le   c_2 U(\d(x)).
\end{equation}
Therefore, by using weak scaling of $\phi$ in the penultimate asymptotic equality,
\begin{eqnarray*}
I_1 &\asymp & \int_{D_1}U(\d(y))\frac{V(|x-y|)^2}{|x-y|^d}\, dy \preceq U(\d(x)) \int_{|y-x|<\d(x)/2} \frac{1}{|x-y|^d\phi(|x-y|^{-2})}\, dy\\
&\preceq  & U(\d(x)) \int_0^{\d(x)}\frac{1}{r\phi(r^{-2})}dr \asymp U(\d(x)) \frac{1}{\phi(\d(x)^{-2})}\asymp U(\d(x))V(\d(x))^2.
\end{eqnarray*}
Finally,  by \eqref{e:wus-U} and the upper weak scaling \eqref{e:wsc-V} of $V$,
\begin{eqnarray*}
\frac{1}{\d(x)}  \int_0^{\d(x)}U(t) V(t)\, dt &\succeq& \frac{U(\d(x))V(\d(x))}{\d(x)}\int_0^{\d(x)}  \frac{V(t)}{V(\d(x))}\, dt\\
&\succeq &\frac{U(\d(x))V(\d(x))}{\d(x)}\int_0^{\d(x)}  \left(\frac{t}{\d(x)}\right)^{\delta_2}\, dt\\
&\asymp & U(\d(x))V(\d(x)).
\end{eqnarray*}

\noindent
{\bf Estimate of $I_2$:} Next, we show that
\begin{equation}\label{e:estimate-I2}
I_2\asymp V(\d(x)) \int_0^{\eta}U(t)V(t)\, dt\, .
\end{equation}
Let $y\in D_2$. Then $r_0 <|y-x|< \mathrm{diam}(D)$ so that $|y-x|\asymp 1$. This implies that $G_D(x,y)\asymp V(\d(x)) V(\d(y))$. Therefore
$$
I_2 \asymp V(\d(x))\int_{D_2} U(\d(y))V(\d(y))\, dy \asymp V(\d(x))\int_{\d(y)<\eta} U(\d(y))V(\d(y))\, dy .
$$
Finally, \eqref{e:estimate-I2} follows by the co-area formula. 

In estimates for $I_3$, $I_4$ and $I_5$ we will use the change of variables formula based on a diffeophormism $\Phi: B(x, r_0)\to B(0,r_0)$ satisfying
\begin{align*}
& \Phi(D\cap B(x,r_0))=B(0,r_0)\cap \{z\in \R^d: \, z\cdot e_d>0\}\\
& \Phi(y)\cdot e_d= \d(y) \ \  \textrm{for any }y\in B(x,r_0), \quad \Phi(x)=\d(x)e_d, 
\end{align*}
see \cite[page 38]{AGCV19}. For the point $z\in \R^d_+=\{z\in \R^d: \, z\cdot e_d>0\}$ we will write $z=(\wt{z}, z_d)$. Several times we also use the following integral:
\begin{equation}\label{e:int}
\int_0^{a}\frac{s^{d-2}}{(b+s)^d}\, ds= \frac{(1+b/a)^{1-d}}{b(d-1)}\, ,\quad a,b>0.
\end{equation}

\noindent 
{\bf Estimate of $I_3$:} It holds that  
\begin{equation}\label{e:estimate-I3}
I_3\asymp \frac{V(\d(x))}{\d(x)}\int_0^{\d(x)}U(t)V(t)\, dt\, .
\end{equation}
To see this, take $y\in D_3$. Then $\d(y)\le \d(x)/2$ implying $|x-y|\ge \d(x)/2$, and thus
\begin{equation}\label{e:G-I3}
G_D(x,y)\asymp \frac{V(\d(x))}{V(|x-y|)}\frac{V(\d(y))}{V(|x-y|)}\frac{V(|x-y|)^2}{|x-y|^d}=\frac{V(\d(x))V(\d(y))}{|x-y|^d}\, .
\end{equation}
Therefore
\begin{eqnarray*}
I_3 &\asymp & V(\d(x))\int_{D_3}\frac{U(\d(y))V(\d(y))}{|x-y|^d}\, dy\\
&\asymp& V(\d(x)) \int_{\{0<z_d<\d(x)/2 \}\cap B(0, r_0)} \frac{U(z_d)V(z_d)}{(|\d(x)-z_d|+|\wt{z}|)^d}\, dz\\
&\asymp & V(\d(x))\int_{|\wt{z}|<r_0} \int_0^{\d(x)/2} \frac{U(z_d)V(z_d)}{(|\d(x)-z_d|+|\wt{z}|)^d}\, dz_d\,  d\wt{z}\\
&\asymp & V(\d(x)) \int_0^{r_0} t^{d-2} \int_0^{\d(x)/2}\frac{U(z_d)V(z_d)}{(|\d(x)-z_d|+t)^d}\, dz_d\, dt\\
&=& V(\d(x)) \int_0^{r_0/\d(x)} s^{d-2} \int_0^{1/2} \frac{U(\d(x)h)V(\d(x)h)}{\big((1-h)+s\big)^d}\, dh\, ds \\
&\asymp & V(\d(x)) \int_0^{r_0/\d(x)} \frac{s^{d-2}}{(1+s)^d}\, ds \int_0^{1/2} U(\d(x)h)V(\d(x)h)\, dh \\ 
&=&V(\d(x)) \frac{\left(1+\d(x)/r_0\right)^{1-d}}{d-1} \int_0^{1/2} U(\d(x)h)V(\d(x)h)\, dh \\
&\asymp & V(\d(x)) \int_0^{1/2} U(\d(x)h)V(\d(x)h)\, dh \\
&=& \frac{V(\d(x))}{\d(x)}\int_0^{\d(x)/2}U(t) V(t)\, dt. 
\end{eqnarray*}
This proves the upper bound in \eqref{e:estimate-I3}. For the lower bound, note that by  the upper weak scaling \eqref{e:wsc-V} of $V$ and the almost nonincreasing condition \eqref{e:wus-U}, we have 
\begin{eqnarray*}
\int_0^{\d(x)/2}U(t) V(t)\, dt&=&2\int_0^{\d(x)} U(t/2)V(t/2)\, dt \ge 2 \int_0^{\d(x)}  c_3 U(t)   \wt a_1^{-1}2^{-\delta_1} V(t)\, dt\\
&=& c_4 \int_0^{\d(x)}U(t)V(t)\, dt\, .
\end{eqnarray*}

\noindent
{\bf Estimate of $I_4$:} By applying the same change of variables as in $I_3$, we show that 
\begin{equation}\label{e:estimate-I4}
I_4\asymp V(\d(x))\int_{3\d(x)/2}^{\eta} \frac{U(t)V(t)}{t}\, dt\, .
\end{equation}
Let $y\in D_4$. Then  $|x-y|\ge \d(x)/2$ and $|x-y|\ge \d(y)/3$, hence $G_D(x,y)$ is of the form \eqref{e:G-I3}. By following the first five lines in the computation of $I_3$, we arrive at
\begin{eqnarray*}
I_4 &\asymp &V(\d(x))\int_0^{r_0/\d(x)}s^{d-2} \int_{3/2}^{\eta/\d(x)} \frac{U(\d(x)h)V(\d(x)h)}{\big((h-1)+s\big)^d}\, dh \, ds\\
&\asymp & V(\d(x)) \int_{3/2}^{\eta/\d(x)}\frac{U(\d(x)h)V(\d(x)h)}{h-1} \int_0^{\frac{r_0}{(h-1)\d(x)}}\frac{r^{d-2}}{(1+r)^d}\, dr \, dh\\
&=& V(\d(x)) \int_{3/2}^{\eta/\d(x)}\frac{U(\d(x)h)V(\d(x)h)}{h-1}  \frac{\left(1+(h-1)\d(x)/r_0\right)^{1-d}}{d-1}\, dh\\
&\asymp& V(\d(x)) \int_{3/2}^{\eta/\d(x)}\frac{U(\d(x)h)V(\d(x)h)}{h-1} \, dh\\
&\asymp & V(\d(x)) \int_{3/2}^{\eta/\d(x)}\frac{U(\d(x)h)V(\d(x)h)}{h} \, dh\\
&=& V(\d(x)) \int_{3\d(x)/2}^{\eta}\frac{U(t)V(t)}{t}\, dt\, .
\end{eqnarray*}

\noindent
{\bf Estimate of $I_5$:}  Under  the almost nonincreasing condition \eqref{e:wus-U} and the doubling condition \eqref{e:doubling-U} it holds that
\begin{equation}\label{e:estimate-I5}
I_5 \preceq U(\delta_D(x))V(\delta_D(x))^2 \preceq \frac{V(\d(x))}{\d(x)}\int_0^{\d(x)}U(t)V(t)\, dt\, .
\end{equation}
Indeed, let $y\in D_5$. Then $|x-y|>\d(x)/2>\d(y)/3$, hence  $G_D(x,y)$ is of the form \eqref{e:G-I3}. Also, the estimate \eqref{e:estimate-I1b} and the analogous one with $V$  hold true. Therefore
\begin{eqnarray*}
I_5 &\asymp & V(\d(x))\int_{D_5} \frac{U(\d(y))V(\d(y))}{|x-y|^d}\, dy\\
&\preceq & U(\d(x))V(\d(x))^2 \int_{D_5}\frac{1}{|x-y|^d}\, dy\, .
\end{eqnarray*}
It is shown in \cite[page 42]{AGCV19} that the last integral is comparable to 1. This proves the first approximate inequality in \eqref{e:estimate-I5}, while the second was already proved in the estimate of $I_1$. 

The proof is finished by noting that $I_1+I_5\preceq I_3$.  \qed

\begin{lemma}\label{l:away-from-bdry}
Let $\eta<\epsilon$ and assume that  conditions \textbf{(U)} hold true. There exists $c(\eta)>0$ such that for any $x\in D$ satisfying $\delta_D(x)\ge\eta/2$, 
\begin{equation}\label{e:away-from-bdry}
G_D\big(U(\delta_D)\1_{(\delta_D<\eta)}\big)(x)\le c(\eta)\, .
\end{equation}
\end{lemma}
\proof Fix $x\in D$ as in the statement and define
\begin{eqnarray*}
D_1&=&\{y:\, \d(y)<\eta/4\},\\
D_2&=&\{y:\, \eta/4 \le \d(y) <\eta\}.
\end{eqnarray*}
Then
$$
G_D\big(U(\delta_D)\1_{(\delta_D<\eta)}\big)(x)=\sum_{j=1}^2 \int_{D_j}G_D(x,y)U(\delta_D(y))\, dy =:\sum_{j=1}^2 J_j.
$$

\noindent 
{\bf Estimate of $J_1$:} 
We show that 
\begin{equation}\label{e:estimate-J1}
J_1\preceq \frac{1}{\eta} \int_0^{\eta}U(t)V(t)\, dt.
\end{equation}
Let $y\in D_1$. Then $\d(y)<\eta/4\le\d(x)/2$, hence by using $|x-y|\ge \d(x)-\d(y)$ we have that $|x-y|>\d(y)$ and $|x-y|>\d(x)/2$. This implies that $G_D(x,y)$ satisfies \eqref{e:G-I3}. Therefore,
$$
J_1\asymp V(\d(x))\int_{D_1}\frac{U(\d(y))V(\d(y))}{|x-y|^d}\, dy.
$$
By using the co-area formula we get (below $dy$ denotes the Hausdorff measure on $\{\d(y)=t\}$)
$$
J_1\asymp \int_0^{\eta/4}U(t)V(t)\left(\int_{\d(y)=t}\frac{1}{|x-y|^d}\, dy\right) dt. 
$$
The inner integral is estimated as follows: For $\d(y)=t$ it holds that $|x-y|\ge \d( x)-t$, hence $|x-y|^{-d}\le (\d(x)-t)^{-d}$. The Hausdorff measure of $\{\d(y)=t\}$ is larger than or equal to the Hausdorff measure of the sphere around $x$ of radius $\d(x)-t$ which is comparable to $(\d(x)-t)^{d-1}$. This implies that the inner integral is estimated from above by a constant times $(\d(x)-t)^{-1}$. Thus
$$
J_1\preceq \int_0^{\eta/4}U(t)V(t) (\d(x)-t)^{-1}\, dt.
$$
If $t<\eta/4$, then $t<\d(x)/2$, implying $\d(x)/2<\d(x)-t<\d(x)$. Therefore,
$$
J_1\preceq \frac{1}{\d(x)}\int_0^{\eta/4}U(t)V(t)\, dt \preceq \frac 2\eta \int_0^{\eta}U(t)V(t)\, dt\, .
$$

\noindent
{\bf Estimate of $J_2$:}
It holds that 
\begin{equation}\label{e:estimate-J2}
J_2\preceq U(\eta/4).
\end{equation}
Let $y\in D_2$. By the almost nonincreasing condition \eqref{e:wus-U} we have $U(\d(y))\le c_1 U(\eta/4)$, hence
\begin{eqnarray*}
J_2 &\preceq & \int_{\eta/4< \d(y) <\eta}U(\d(y))\frac{V(|x-y|)^2}{|x-y|^d}\, dy \preceq U(\eta/4)\int_{\eta/4< \d(y) <\eta}\frac{V(|x-y|)^2}{|x-y|^d}\, dy\\
&\le & U(\eta/4)\int_{B(x, 2\mathrm{diam}(D))} \frac{V(|x-y|)^2}{|x-y|^d}\, dy \preceq U(\eta/4).
\end{eqnarray*}
The last estimate uses the fact that the integral is not singular.

By putting together estimates for $J_1$ and $J_2$, we see that there exists $c_2>0$ such that
$$
G_D\big(U(\delta_D)\1_{(\delta_D<\eta)}\big)(x)\le c_2\left(\frac{1}{\eta} \int_0^{\eta}U(t)V(t)\, dt +U(\eta/4)\right)=: c(\eta).
$$
\qed

\noindent
\textbf{Proof of Proposition \ref{p:green-potential-estimate}: } First we prove the statement under conditions \textbf{(U)}. Fix some $\eta<\epsilon$ and treat it as a constant. Note that on $\{\d(y)\ge \eta\}$ it holds that $U$ is bounded (by the assumption  (U4)). Therefore
\begin{equation}\label{e:gpe-b}
G_D(U(\d)\1_{(\d\ge \eta)})(x)\asymp G_D \1(x)\asymp V(\d(x))\, .
\end{equation}
By Lemma \ref{l:away-from-bdry}, if $\d(x)\ge \eta/2$, then $G_D\big(U(\delta_D)\1_{(\delta_D<\eta)}\big)(x)\le c(\eta)$. Hence, 
$$
G_D(U(\d))(x)\asymp 1, \quad \d(x)\ge \eta/2.
$$
Since for $\d(x)\ge \eta/2$ the right-hand side in \eqref{e:green-potential-estimate} is also comparable to 1, this proves the claim for this case. Assume now that $\d(x)<\eta/2$. By Lemma \ref{l:close-to-bdry} and \eqref{e:gpe-b} we have that
\begin{eqnarray*}
G_D(U(\delta_D))(x)&=&G_D(U(\d)\1_{(\d< \eta)})(x)+G_D(U(\d)\1_{(\d\ge \eta)})(x)\\
&\asymp &\frac{V(\d(x))}{\d(x)}\int_0^{\d(x)}U(t)V(t)\, dt + V(\d(x))\int_0^{\eta}U(t)V(t)\, dt \nonumber \\
& &+\ V(\d(x))\int_{3\d(x)/2}^{\eta} \frac{U(t)V(t)}{t}\, dt\ + V(\d(x))\\
&\asymp & \frac{V(\d(x))}{\d(x)}\int_0^{\d(x)}U(t)V(t)\, dt + V(\d(x))\int_{\d(x)}^{\eta} \frac{U(t)V(t)}{t}\, dt\, .
\end{eqnarray*}
Clearly, in the last integral we can replace $\eta$ by $\mathrm{diam}(D)$.

Lastly, assume that the function $U$ is bounded on every bounded subset of $(0,\infty)$. Obviously, by \eqref{e:gpe-b},
\begin{align*}
G_D(U(\delta_D))(x)\preceq G_D\1(x)\asymp V(\d(x)).
\end{align*}
On the other hand, analogously as in \eqref{e:estimate-I2},
\begin{align*}
G_D(U(\delta_D))(x)\geq \int_{D_2}U(\d(y))G_D(x,y)\,dy\asymp V(\d(x))\int_0^{\eta}U(t)V(t)\, dt\, .
\end{align*}

\qed

\noindent
\textbf{Proof of Proposition \ref{p:poisson-potential-estimate}:} Fix $\eta<\epsilon$. Let $x\in D$ and $r_0>\d(x)+\eta$. We split $D^c$ into three parts,
\begin{align*}
D_1&=\{z\in D^c:\,\dc(z)\geq \eta\}\\
D_2&=\{z\in D^c\cap B(x,r_0):\,\dc(z)<\eta\}\\
D_3&=\{z\in D^c\setminus B(x,r_0):\,\dc(z)< \eta\}
\end{align*} 
and apply \eqref{e:g_G} to get that
\begin{align*}
P_Dg(x)&\asymp \int_{D_1}\widetilde U(\dc(z))P_D(x,z)dz+\int_{D_2}\widetilde U(\dc(z))P_D(x,z)dz+\int_{D_3}\widetilde U(\dc(z))P_D(x,z)dz\\
&=:I_1+I_2+I_3.
\end{align*}

\noindent
{\bf Estimate of $I_1$:}
For $z\in D^c$ such that $\dc(z)\geq \eta$, the estimate \eqref{e:P_D-estimate} is equivalent to
\begin{equation*}
P_D(x,z)\asymp \frac{V(\d(x))}{V(\dc(z))^2\dc(z)^d}.
\end{equation*}
By applying this estimate and the co-area formula to $I_1$, we arrive to
\begin{align*}
I_1&\asymp V(\d(x))\int_{D_1}\frac{\widetilde U(\dc(z))}{V(\dc(z))^2\dc(z)^d}dz\\
&\asymp V(\d(x))\int_\eta^\infty\frac{\widetilde U(t)}{V(t)^2t^d}\int_{D^c} \1_{\delta_D(w)=t}\,  dw dt\\
&\asymp V(\d(x))\int_\eta^\infty\frac{\widetilde U(t)}{V(t)^2t} dt.
\end{align*}
As before, $dw$ in the first two lines denotes the Hausdorff measure on $\d(w)=t$ and we used that 
$$
|\{w\in D^c:\, \dc(w)=t\}|\asymp t^{d-1},\ t\geq\eta.  
$$

\noindent
{\bf Estimate of $I_2$:}
First note that for $z\in D^c\cap B(x,r_0)$ estimate \eqref{e:P_D-estimate} implies that
$$
P_D(x,z)\asymp \frac{V(\d(x))}{V(\dc(z))|x-z|^d}.
$$
Next, as in the proof of Lemma \ref{l:close-to-bdry} we will use the change of variables formula based on a diffeophormism $\Phi: B(x, r_0)\to B(0,r_0)$ satisfying
\begin{align*}
& \Phi({\overline D}^c\cap B(x,r_0))=B(0,r_0)\cap \{w\in \R^d: \, w\cdot e_d<0\}\\
& |\Phi(z)\cdot e_d|= \d(z) \ \  \textrm{for any }z\in {\overline D}^c\cap B(x,r_0), \quad \Phi(x)=\d(x)e_d\, .
\end{align*}
Similarly as before, for the point $w\in \R^d_-=\{w\in \R^d: \, w\cdot e_d<0\}$ we will write $w=(\wt{w}, w_d)$. Therefore, by the change of variables given by the diffeomorphism $\Phi$ it follows that
\begin{align*}
I_2&\asymp V(\d(x))\int_{D_2}\frac{\widetilde U(\dc(z))}{V(\dc(z))|x-z|^d}dz\\
&\asymp V(\d(x))\int_{\{w\in B(0,r_0):-\eta<w_d<0\}}\frac{\widetilde U(-w_d)}{V(-w_d)(|\d(x)-w_d|+|\widetilde{w}|)^d}dw.
\end{align*}
Next, we apply the substitution $w_d=-t$ and switch to polar coordinates for $\widetilde{w}$ to obtain that   
\begin{align*}
I_2&\asymp V(\d(x))\int_0^\eta\frac{\widetilde U(t)}{V(t)}\int_0^{r_0}\frac{s^{d-2}}{(\d(x)+t+s)^d}ds\,dt\\
&\overset{\eqref{e:int}}{\asymp} V(\d(x))\int_0^\eta\frac{\widetilde U(t)}{V(t)(\d(x)+t)}dt\\
&\leq \frac{V(\d(x))}{\d(x)}\int_0^\eta\frac{\widetilde U(t)}{V(t)}dt.
\end{align*}
\noindent
{\bf Estimate of $I_3$:}
Lastly, note that for $z\in D^c\setminus B(x,r_0)$ such that $\d(z)<\eta$, estimate \eqref{e:P_D-estimate} is equivalent to
$$
P_D(x,z)\asymp \frac{V(\d(x))}{V(\dc(z))}.
$$
Therefore, similarly as in the estimate of $I_1$ we have
\begin{align*}
I_3&\asymp V(\d(x))\int_{D_3}\frac{\widetilde U(\dc(z))}{V(\dc(z))}dz\\
&\asymp V(\d(x))\int_{0}^\eta\frac{\widetilde U(t)}{V(t)}\int_{ D^c\setminus B(x,r_0)}\1_{\d(w)=t}\, dw dt\\
&\asymp V(\d(x))\int_{0}^\eta\frac{\widetilde U(t)}{V(t)}dt.
\end{align*}
Since for $t<\eta$ we have that $\delta_D(x)+t<\mathrm{diam}(D)+\eta$, it follows that $I_3\preceq I_2$.

This proves that
$$
P_Dg(x)\asymp V(\delta_D(x))\left(\int_0^\eta\frac{\widetilde U(t)}{V(t)(\d(x)+t)}dt+\int_\eta^\infty \frac{\widetilde U(t)}{V(t)^2t}dt\right),\ x\in D.
$$
By fixing $\eta$ and noting that 
$$
\int_{\eta}^{\mathrm{diam}(D)}\frac{\widetilde U(t)}{V(t)(\d(x)+t)}dt+ \int_\eta^\infty \frac{\widetilde U(t)}{V(t)^2t}dt\asymp 1
$$
we  obtain \eqref{e:P_Dg}. Inequality \eqref{e:P_Dg2} follows immediately.
\qed

\bigskip

{\bf Ivan Bio\v{c}i\'c}

Department of Mathematics, Faculty of Science, University of Zagreb, Zagreb, Croatia,

Email: \texttt{ibiocic@math.hr}

\bigskip

{\bf Zoran Vondra\v{c}ek}

Department of Mathematics, Faculty of Science, University of Zagreb, Zagreb, Croatia,

Email: \texttt{vondra@math.hr}

\bigskip
{\bf Vanja Wagner}

Department of Mathematics, Faculty of Science, University of Zagreb, Zagreb, Croatia,

Email: \texttt{wagner@math.hr}
\end{document}